\documentclass[12pt,a4paper]{article}
\usepackage{amssymb,amsmath}
\usepackage{amsfonts}
\usepackage{a4}
\usepackage{tikz}
\usepackage{verbatim}
\usepackage{caption}
\usepackage{float}
\newtheorem{prop}{Proposition}
\newtheorem{theorem}{Theorem}
\newtheorem{lemma}{Lemma}
\newtheorem{cor}{Corollary}
\newtheorem{remark}{Remark}

\begin{document}
\baselineskip 18pt \title{\bf \ Character triples and Shoda pairs}
\author{ Gurmeet K. Bakshi and Gurleen Kaur{\footnote {Research supported by NBHM, India, is gratefully acknowledged} \footnote{Corresponding author}} \\ {\em \small Centre for Advanced Study in
Mathematics,}\\
{\em \small Panjab University, Chandigarh 160014, India}\\{\em
\small email: gkbakshi@pu.ac.in and gurleenkaur992gk@gmail.com  } }
\date{}
{\maketitle}
\begin{abstract}\noindent
{In this paper, a construction of Shoda pairs using character triples is given for a large class of monomial groups including abelian-by-supersolvable and subnormally monomial groups. The computation of primitive central idempotents and the structure of simple components of the rational group algebra for groups in this class are also discussed. The theory is illustrated with examples. }
\end{abstract}\vspace{.25cm}
 \vspace{.25cm} 
\section{Introduction}
  \indent Given a finite group $G$,  Shoda (\cite{CR}, Corollary 45.4) gave a criterion to determine whether an induced monomial representation of  $G$ is irreducible or not. Olivieri, del R{\'{\i}}o and Sim{\'o}n  \cite{OdRS04} rephrased Shoda's theorem as follows: \begin{quote} If $\chi$ is a linear character of a subgroup $H$ of $G$ with kernel $K$, then the induced character $\chi^{G}$ is irreducible if, and only if, the following hold:\vspace{.2cm}\\(i) $K\unlhd H$, $H/K$ is cyclic;\\(ii) if $g \in G$ and $[H,g]\cap H\subseteq K$, then $g \in H$.\end{quote} A pair $(H,K)$ of subgroups of $G$ satisfying (i) and (ii) above is called a Shoda pair of $G$. For $K\unlhd H\leqslant G$, define:$$\widehat{H}:=\frac{1}{|H|}\displaystyle\sum_{h \in H}h,$$ $$\varepsilon(H,K):=\left\{\begin{array}{ll}\widehat{K}, & \hbox{H=K;} \\\prod(\widehat{K}-\widehat{L}), & \hbox{otherwise,}\end{array}\right.$$ where $L$ runs over all the minimal normal subgroups of $H$ containing $K$ properly, and $$e(G,H,K):= {\rm~the~sum~of~all~the~distinct~}G{\rm {\tiny{\operatorname{-}}} conjugates~of~}\varepsilon(H,K).$$ An important feature (\cite{OdRS04}, Theorem 2.1) of a Shoda pair $(H,K)$ of $G$ is that there is a rational number  $\alpha$, necessarily unique, such that $\alpha e(G,H,K)$ is a primitive central idempotent of the rational group algebra $\mathbb{Q}G$, called the primitive central idempotent of $\mathbb{Q}G$ realized  by the Shoda pair $(H,K)$. We'll denote this $\alpha$ by $\alpha_{(G, H, K)}$. For monomial groups, all the primitive central idempotents of $\mathbb{Q}G$ are realized by Shoda pairs of $G$. For the Shoda pair $(H,K)$ of $G$, the case when $e(G,H,K)$ is a primitive central idempotent of $\mathbb{Q}G$, is of special interest, thus, leading to the following definition of a strong Shoda pair. A strong Shoda pair \cite{OdRS04} of $G$ is a pair $(H, K)$ of subgroups of $G$ satisfying the following conditions: \begin{quote}(i) $K\unlhd H \unlhd N_{G}(K)$;\\(ii) $H/K$ is cyclic and a maximal abelian subgroup of $N_{G}(K)/K$;\\(iii) the distinct $G$-conjugates of
$\varepsilon(H,K)$ are mutually orthogonal.\end{quote} In \cite{OdRS04}, it is proved that if $(H,K)$ is a strong Shoda pair of $G$, then it is also a Shoda pair of $G$ and $e(G,H,K)$ is a primitive central idempotent of $\mathbb{Q}G.$ The groups $G$ such that all the primitive central idempotents of $\mathbb{Q}G$ are realized by strong Shoda pairs of $G$ are termed as strongly monomial groups. Examples of such groups include abelian-by-supersolvable groups (\cite{OdRS04}, Theorem 4.4). The main reason for defining strong Shoda pairs in \cite{OdRS04} was that the authors were able to provide a description of the structure of the simple component $\mathbb{Q}Ge(G,H,K) $ of $\mathbb{Q}G$ for a strong Shoda pair $(H,K)$ of $G$. \par The work in \cite{OdRS04} thus leads to the problem of computing Shoda pairs of a given finite group $G$ and to provide a description of the structure of the simple components of $\mathbb{Q}G$ corresponding to the primitive central idempotents realized by them. The interest is in fact in providing a method to obtain a set $\mathcal{S}$ of Shoda pairs of $G$ such that the mapping $(H, K)$ $\mapsto$ $\alpha_{(G, H, K)}e(G,H,K)$ defines a bijection from $\mathcal{S}$ to the set of all primitive central idempotents of $\mathbb{Q}G$ realized by Shoda pairs of $G$. Such a set $\mathcal{S}$ is called a complete and irredundant set of Shoda pairs  of $G$, and has recently been provided by the first author with Maheshwary \cite{BM} for normally monomial groups. For the work in this direction, also see \cite{BKP} and \cite{BM1}.\par In this paper, we plan to study the problem for the class $\mathcal{C}$ of all finite groups $G$ such that all the subgroups and quotient groups of $G$ satisfy the following property: either it is abelian or it contains a non central abelian normal subgroup. The groups in $\mathcal{C}$ are known to be monomial (\cite{BH}, Lemma 24.2). However, we have noticed that $\mathcal{C}$ is not contained in the class of strongly monomial groups. Huppert (\cite{BH}, Theorem 24.3) proved that $\mathcal{C}$ contains all the groups $G$ for which there is a solvable normal subgroup $N$ of $G$ such that all Sylow subgroups of $N$ are abelian and $G/N$ is supersolvable. In particular, $\mathcal{C}$ contains abelian-by-supersolvable groups. In view of an important criterion of subnormally monomial groups given in \cite{HowI} and \cite{HowII}, we have shown in section 2 that $\mathcal{C}$ also contains all subnormally monomial groups and, in particular, normally monomial groups. Our aim is to extend the work to the class $\mathcal{C}$.\par  An important tool which has turned out to be useful is Isaacs's notion of character triples together with Clifford's correspondence theorem. Following Isaacs (\cite{IM}, p.186), we have defined $N$-linear character triples of $G$ for a normal subgroup $N$ of $G$. In view of Clifford's correspondence theorem (\cite{IM}, Theorem 6.11), for each $N$-linear character triple of $G$, we have defined its direct Clifford correspondents, which is another set of $N$-linear character triples of $G$ with useful properties proved in Theorem \ref{t1} of section 3. With its help, we have given, in section 4, a construction of Shoda pairs of groups in $\mathcal{C}$. For each normal subgroup $N$ of $G$, we have constructed a rooted directed tree $\mathcal{G}_{N}$, whose particular leaves correspond to Shoda pairs of $G$, if $G \in \mathcal{C}$ (Theorem \ref{t2}). We have also explored the condition for the collection of Shoda pairs corresponding to these leaves of $\mathcal{G}_{N}$ as $N$ runs over all the normal subgroups of $G$ to be complete and irredundant. In section 5, we have given a new character free expression of $\alpha_{(G,H,K)}$, where $(H,K)$ is a Shoda pair of $G$ corresponding to a leaf of $\mathcal{G}_{N}$. This expression is in terms of the directed path from the root to the corresponding leaf and enables us to provide a necessary and sufficient condition for $e(G,H,K)$ to be a primitive central idempotent of $\mathbb{Q}G$. In section 6, we generalize Proposition 3.4 of \cite{OdRS04} and determine the structure of the simple components of $\mathbb{Q}G$ for $G \in \mathcal{C}$. Finally, in section 7, we provide illustrative examples.

\section{The class $\mathcal{C}$ of monomial groups} Throughout this paper, $G$ denotes a finite group. By $H\leq G$, $H \lneq G$, $H\unlhd G$, we mean, respectively, that $H$ is a subgroup, proper subgroup, normal subgroup of $G$. Denote by $[G:H]$, the index of $H$ in $G$. Also $N_{G}(H)$ denotes the normalizer of $H$ in $G$ and $\operatorname{core}_{G}(H)=\bigcap_{x \in G}xHx^{-1}$ is the largest normal subgroup of $G$ contained in $H$. For $x,~y \in G$, $[x,y]=x^{-1}y^{-1}xy$ is the commutator of $x$ and $y$, and $\operatorname{Cen}_{G}(x)=\{g \in G~|~gx=xg\}$ is the centralizer of $x$ in $G$. Denote by $\operatorname{Irr}G$, the set of all complex irreducible characters of $G$. For a character $\chi$ of $G$, $\operatorname{ker}\chi =\{ g \in G~|~\chi(g)=\chi(1)\}$ and $\mathbb{Q}(\chi)$ denotes the field obtained by adjoining to $\mathbb{Q}$ the character values $\chi(g)$, $g \in G$. If $\psi$ is a character of a subgroup $H$ of $G$ and $x \in G$, then $\psi^{x}$ is the character of $H^{x}=x^{-1}Hx$ given by $\psi^{x}(g)=\psi(xgx^{-1})$, $g \in H^{x}$. Denote by $\psi^{G}$, the character $\psi$ induced to $G$. For a subgroup $A$ of $H$, $\psi_{A}$ denotes the restriction of $\psi$ to $A$. \par Let $\mathcal{C}$ denote the class of all finite groups $G$ such that all the subgroups and quotient groups of $G$ satisfy the following property: either it is abelian or it contains a non central abelian normal subgroup. It follows from Lemma 24.2 of \cite{BH} that the groups in $\mathcal{C}$ are monomial. Recall that a finite group is monomial if every complex irreducible character of the group is induced by a linear character of a subgroup. In this section, we compare $\mathcal{C}$ with the following classes of groups:
 \vspace{.3cm}\\ $\begin{array}{lll} \operatorname{Ab} &:& {\rm Abelian~groups}\vspace{.2cm}\\
  \operatorname{Nil}&:& {\rm Nilpotent~groups}\vspace{.2cm}\\
  \operatorname{Sup}&:& {\rm Supersolvable~groups}\vspace{.2cm}\\
 \mathcal{A}&:&{\rm Solvable~groups~}{\rm with~all~the~Sylow~subgroups~abelian}\vspace{.2cm}\\ \operatorname{nM} &:& {\rm Normally ~ monomial~groups,~i.e.,~groups~with~all~the~complex}
\vspace{.1cm}\\&&{\rm irreducible ~characters}{\rm~induced~from~linear~characters~of~normal}\vspace{.1cm}\\&&{\rm subgroups}\vspace{.2cm}\\ \operatorname{sM}&:&{\rm Subnormally~monomial~groups,~i.e.,~groups~with~all~the~complex}
\vspace{.1cm}\\&&{\rm irreducible~characters~induced~from~linear~characters~of~subnormal}\vspace{.1cm}\\ &&{\rm subgroups}\vspace{.2cm}\\\operatorname{stM}&:&{\rm Strongly~monomial~groups}\vspace{.2cm}\\ X&:&{\rm Solvable~groups~}G{\rm~satisfying~the~following~condition}:{\rm~For~all~}\vspace{.1cm}\\&&{\rm primes~}p{\rm~dividing~the~order~of~}G{\rm~and~for~all~subgroups}
~A{\rm~of~}G,\vspace{.1cm}\\&& O^{p}(A),~{\rm the~unique~smallest~normal~subgroup~of~}A{\rm~such~that~}\vspace{.1cm}\\&& A/O^{p}(A)~{\rm is~a~}p\tiny{\operatorname{-}}{\rm group}, {\rm~has~no~} {\rm central~} p\tiny{\operatorname{-}}{\rm factor}\vspace{.2cm}\\\mathcal{X}\tiny{\operatorname{-}}\operatorname{by}\tiny{\operatorname{-}}\mathcal{Y}&:&{\rm Groups~}G{\rm~such~that~there~exist~ a~normal~subgroup~} N{\rm~with~}N\in \mathcal{X}\vspace{.1cm}\\&& {\rm and~} G/N \in \mathcal{Y}. \end{array}$
\vspace{.3cm}\\We prove the following:\begin{prop}\label{p0} The following statements hold:\begin{description}\item[(i)]$\operatorname{Ab-by-Nil}$ $\subseteq$ $\operatorname{Ab-by-Sup}$ $\subseteq$ $\mathcal{A}$-$\operatorname{by-Sup}$ $\subseteq$ $\mathcal{C}$;\item[(ii)] $\operatorname{(}\operatorname{nM}$ $\cup$ $\operatorname{Ab-by-Nil}\operatorname{)}$ $\subseteq$ $\operatorname{sM}$  $\subseteq$ X $\subseteq$ $\mathcal{C}$; \item[(iii)] $\mathcal{A}$-$\operatorname{by-Sup}$ $\nsubseteq$ X;  \item[(iv)] X $\nsubseteq$ $\mathcal{A}$-$\operatorname{by-Sup}$; \item[(v)]  $\mathcal{C}$ $\nsubseteq$ $\operatorname{stM}$. \end{description}\end{prop}
{\bf Proof.} (i) Clearly, Ab-by-Nil $\subseteq$ Ab-by-Sup $\subseteq$ $\mathcal{A}$-by-Sup. From (\cite{BH}, Theorem 24.3), it follows that $\mathcal{A}$-by-Sup $\subseteq \mathcal{C}$. This proves (i).\vspace{.1cm}\\ (ii) It is obvious that nM $\subseteq$ sM. We now show that Ab-by-Nil $\subseteq$ sM. Let $G \in$ Ab-by-Nil. Let $A$ be a normal abelian subgroup of $G$ such that $G/A$ is nilpotent. Let $\chi \in \operatorname{Irr} G$. It is already known that $\chi$ is monomial. By (\cite{BH}, Lemma 24.8), there exists a subgroup $H$ of $G$ containing $A$ such that $\chi$ is induced from a linear character on $H$. As $H/A$ is a subgroup of the nilpotent group $G/A$, it is subnormal in $G/A$. Consequently, $H$ is subnormal in $G$. This proves that $G \in$ sM. Next, by (\cite{HowII}, Theorem 3.7), we have sM $\subseteq X$. We now show that $X\subseteq \mathcal{C}$. By Lemma 2.6 of \cite{HowI}, $X$ is closed under taking subgroups and factor groups. Thus to prove that $X\subseteq \mathcal{C}$, we only need to show that every non abelian group in $X$ contains a non central abelian normal subgroup. Let $G \in X.$ If $G$ is nilpotent, then clearly it has the desired property. If $G$ is not nilpotent, then Lemma 2.7 of \cite{HowI} implies that $\sigma(G)$, the socle of $G$, is non central in $G$. Also, in view of (\cite{IM1}, Lemma 3.11, Problem 2A.5), $\sigma(G)$ is abelian, as $G$ is solvable. Hence $G$ has the desired property and it follows that $X\subseteq \mathcal{C}$.\vspace{.1cm}\\(iii) Consider the group $G$ generated by $a, b, c, d$ with defining relations: $a^{2}=b^{3}=c^{3}=d^{3}=1,~a^{-1}ba=b^{-1},~a^{-1}ca=c^{-1},~a^{-1}da=d,~b^{-1}cb=cd,~b^{-1}db=d,~c^{-1}dc=d.$ It is easy to see that $G$ is supersolvable and hence belongs to $\mathcal{A}\tiny\operatorname{-}$by$\tiny\operatorname{-}$Sup. We'll show that $G \not\in X$. Let $CSF$ be the class of all chiefly sub-Frobenius groups, i.e., all finite solvable groups $G$ for which $\operatorname{Cen}_{G}(kL)$ is subnormal in $G$ whenever $kL$ is an element of a chief factor $K/L$ of $G$. It is easy to see that supersolvable groups are chiefly sub-Frobenius, and hence $G \in CSF$. Also it is known (\cite{HowII}, Theorem 3.8) that sM=$CSF\cap X$. Now if $G \in X$, then it follows that $G$ is subnormally monomial. It can be shown that $G$ has an irreducible character of degree $3$, denoted by $\chi$ say. If $G$ is subnormally monomial, then $\chi$ is induced from a linear character of a subnormal subgroup $H$ of $G$. Also, $\chi(1)=3$ implies that $[G:H]=3$, which yields that $H$ is in fact normal in $G$. However, $G$ does not have any normal subgroup of index $3$, a contradiction. This proves that $G \not \in X$ and (iii) follows.\vspace{.1cm}\\(iv) Consider $G$=SmallGroup$(192,1025)$ in GAP\cite{GAP4} library. It can be checked using GAP that $G$ is subnormally monomial and hence belongs to the class $X$ but it does not belong to $\operatorname{\mathcal{A}}$-by-Sup.\vspace{.1cm}\\(v) Simple computations using Wedderga\cite{wedd} reveal that SmallGroup$(1000,86)$ is not strongly monomial. However, it belongs to $\mathcal{C}$. This proves (v). $~\Box$

\section{$N$-linear character triples} Let $G$ be a finite group. Let $H\leq G$ and $(H,A,\vartheta)$ a character triple, i.e., $A\unlhd H$ and $\vartheta \in \operatorname{Irr}A$ is invariant in $H$. We call it to be $N$-linear character triple of $G$, if $\vartheta$ is linear and $\operatorname{ker}\vartheta^{G}=N$. For $N$-linear character triple $(H,A,\vartheta)$ of $G$, denote by ${\operatorname{Irr}}(H|\vartheta)$, the set of all irreducible characters $\psi$ of $H$ which lie above $\vartheta$, i.e., the restriction $\psi_{A}$ of $\psi$ to $A$ has $\vartheta$ as a constituent. Let $\widetilde{\operatorname{Irr}}(H|\vartheta)$ be its subset consisting of those $\psi$ which satisfy $\operatorname{ker}\psi^{G} = N.$ Denote by $\operatorname{\widetilde{Lin}}(H|\vartheta)$, the subset of $\widetilde{\operatorname{Irr}}(H|\vartheta)$ consisting of linear characters. Further, for the character triple $(H,A,\vartheta)$, we fix a normal subgroup $\mathcal{A}_{(H, A, \vartheta)}$ of $H$ of maximal order containing $\operatorname{ker} \vartheta$ such that $\mathcal{A}_{(H, A, \vartheta)}/\operatorname{ker} \vartheta$ is abelian. Note that there may be several choices of such $\mathcal{A}_{(H, A, \vartheta)}$, however, we fix one such choice for a given triple $(H,A,\vartheta)$. Observe that $\mathcal{A}_{(H, A, \vartheta)}/\operatorname{ker} \vartheta$ always contains the center of $H/\operatorname{ker} \vartheta$. However, it can be seen that if $G \in \mathcal{C}$ and $H/\operatorname{ker} \vartheta$ is non abelian, then  $\mathcal{A}_{(H, A, \vartheta)}/\operatorname{ker} \vartheta$ properly contains the center of $H/\operatorname{ker} \vartheta $. We shall later use this observation without any mention. \par Given $N$- linear character triple $(H,A,\vartheta)$ of $G$, we provide a construction of the set $Cl(H,A,\vartheta)$ of another $N$-linear character triples of $G$ required for the purpose of constructing Shoda pairs of $G$.\vspace{.2cm}\\\underline{\bf Construction of $Cl(H,A,\vartheta)$}\vspace{.2cm}\\ Let $\operatorname{\mathcal{A}ut}(\mathbb{C}|\vartheta)$ be the group of automorphisms of the field $\mathbb{C}$ of complex numbers which keep $\mathbb{Q}(\vartheta)$ fixed. For brevity, denote $\mathcal{A}_{(H, A, \vartheta)}$ by $\mathcal{A}$. Consider the action of $\operatorname{\mathcal{A}ut}(\mathbb{C}|\vartheta)$ on $\operatorname{\widetilde{Lin}}(\mathcal{A}|\vartheta)$ by setting $$ ~~~~~~~~~~~~~~~~~\sigma.\varphi= \sigma\circ\varphi,~~~~~~~~~~~~~~~~\sigma\in \operatorname{\mathcal{A}ut}(\mathbb{C}|\vartheta),~\varphi \in \operatorname{\widetilde{Lin}}(\mathcal{A}|\vartheta).$$ Also $H$ acts on $\operatorname{\widetilde{Lin}}(\mathcal{A}|\vartheta)$ by $$ ~~~~~~~~~~~h.\varphi= \varphi^{h},~~~~~~~~~~~~~~~~h\in H,~\varphi \in \operatorname{\widetilde{Lin}}(\mathcal{A}|\vartheta).$$ Notice that the two actions on $\operatorname{\widetilde{Lin}}(\mathcal{A}|\vartheta)$ are compatible in the sense that $$~~~~~~~~~~~~\sigma.(h.\varphi) = h.(\sigma.\varphi),~~~~~~~~~~~~~~~~h\in H,~\sigma \in \operatorname{\mathcal{A}ut}(\mathbb{C}|\vartheta),~\varphi \in \operatorname{\widetilde{Lin}}(\mathcal{A}|\vartheta).$$ This consequently gives an action of $\operatorname{\mathcal{A}ut}(\mathbb{C}|\vartheta)\times H $ on $\operatorname{\widetilde{Lin}}(\mathcal{A}|\vartheta).$ Under this double action, denote by $\operatorname{\mathfrak{Lin}}(\mathcal{A}|\vartheta)$, a set of representatives of distinct orbits of $\operatorname{\widetilde{Lin}}(\mathcal{A}|\vartheta)$. If $H \neq A$, set $$Cl(H,A,\vartheta)=\{(I_{H}(\varphi),\mathcal{A},\varphi)~|~\varphi \in \operatorname{\mathfrak{Lin}}(\mathcal{A}|\vartheta) \},$$ where $I_{H}(\varphi)=\{g \in G~|~\varphi^{g}=\varphi\}$ is the inertia group of $\varphi$ in $H$. For $H=A$, define $Cl(H,A,\vartheta)$ to be an empty set. Note that all the character triples in $Cl(H,A,\vartheta)$ are $N$-linear character triples of $G$ and we call them the direct Clifford correspondents (abbreviated d.c.c.) of $(H,A,\vartheta)$. The name `direct Clifford correspondents' refers to the fact that the characters in $\operatorname{\widetilde{Irr}}(I_{H}(\varphi)|\varphi)$ are Clifford correspondents of the characters in $\operatorname{\widetilde{Irr}}(H|\vartheta)$ in view of the following theorem:
 \begin{theorem}\label{t1} Let $G \in \mathcal{C}$ and $N$ a normal subgroup of $G$. Let $(H,A,\vartheta)$ be $N$-linear character triple of $G$ with $H\neq A$ and $Cl(H,A,\vartheta)$ be as defined above. Let $\mathcal{A} = \mathcal{A}_{(H, A, \vartheta)}$. Then \begin{description} \item[(i)] for any $(I_{H}(\varphi),\mathcal{A},\varphi) \in Cl(H,A,\vartheta),$ the following hold:
  \begin{description}\item[(a)] $A\lneq \mathcal{A}$. Furthermore, $\mathcal{A}=H= I_{H}(\varphi)$ holds, if, and only if, $H/\operatorname{ker} \vartheta$ is abelian; \item[(b)] $I_{H}(\varphi) \unlhd N_{H}(\operatorname{ker}\varphi)$; \item[(c)] the induction $\psi \mapsto \psi^{H}$ defines an injective map from $\operatorname{\widetilde{Irr}}(I_{H}(\varphi)|\varphi)$ to $\operatorname{\widetilde{Irr}}(H|\vartheta)$.\end{description} \item[(ii)] for each $\chi \in \operatorname{\widetilde{Irr}}(H|\vartheta)$, there exists $(I_{H}(\varphi),\mathcal{A},\varphi)\in Cl(H,A,\vartheta)$, $\psi \in \operatorname{\widetilde{Irr}}(I_{H}(\varphi)|\varphi)$ and $\sigma \in \operatorname{\mathcal{A}ut}(\mathbb{C})$ such that $\chi = \sigma \circ \psi^{H}.$  \item[(iii)] if $(I_{H}(\varphi_{1}),\mathcal{A},\varphi_{1}),~(I_{H}(\varphi_{2}),\mathcal{A},\varphi_{2}) \in Cl(H,A,\vartheta)$, $\psi_{1}\in \operatorname{\widetilde{Irr}}(I_{H}(\varphi_{1})|\varphi_{1}),$\linebreak $\psi_{2}\in \operatorname{\widetilde{Irr}}(I_{H}(\varphi_{2})|\varphi_{2})$ and $ \sigma \in \operatorname{\mathcal{A}ut}(\mathbb{C})$ are such that $\psi_{2}^{H}= \sigma \circ \psi_{1}^{H}$, then\linebreak $\varphi_{1}=\varphi_{2} ~(= \varphi~\rm{say})$, and in this case $\psi_{2}=\sigma\circ \psi_{1}^{x}$ for some $x\in N_{H}(\operatorname{ker}\varphi)$.  \end{description} \end{theorem}{\bf Proof.} (i) We first show that \begin{equation}\label{e1} A\lneq \mathcal{A}. \end{equation}
  If $H/\operatorname{ker}\vartheta$ is abelian, then clearly $\mathcal{A} = H$ and therefore the above equation holds trivially, as $A\neq H$. If $H/\operatorname{ker}\vartheta$ is non abelian, then $\mathcal{A}/\operatorname{ker}\vartheta$ properly contains the centre of $H/\operatorname{ker}\vartheta$. However, $\vartheta$ being invariant in $H$, it follows that $A/\operatorname{ker}\vartheta$ is contained in the centre of $H/\operatorname{ker}\vartheta$. Therefore, eqn (\ref{e1}) follows. This proves (a). Next, consider $\alpha \in I_{H}(\varphi)$ and $\beta \in N_{H}(\operatorname{ker} \varphi).$ Then $$\beta^{-1}\alpha\beta \in I_{H}(\varphi)$$ if, and only if, \begin{equation}\label{e2}[\beta^{-1}\alpha\beta, x]=\beta^{-1}[\alpha,\beta x \beta^{-1}]\beta \in \operatorname{ker} \varphi~{\rm for~all}~x \in \mathcal{A}.\end{equation} However, if $x \in \mathcal{A}$, then $\beta x \beta^{-1} \in \mathcal{A}$, as $\mathcal{A} \unlhd H$. This gives $[\alpha,\beta x \beta^{-1}] \in \operatorname{ker}\varphi$, as $\alpha \in I_{H}(\varphi)$. Consequently, $\beta$ being in $N_{H}(\operatorname{ker} \varphi)$, eqn (\ref{e2}) follows. This proves (b). In view of Clifford's correspondence theorem (\cite{IM}, Theorem 6.11), $\psi\mapsto \psi^{H}$ defines an injective map from $\operatorname{{Irr}}(I_{H}(\varphi)|\varphi)$ to $\operatorname{{Irr}}(H|\vartheta)$. It can be checked that, under this map, $\operatorname{\widetilde{Irr}}(I_{H}(\varphi)|\varphi)$ is mapped to $\operatorname{\widetilde{Irr}}(H|\vartheta)$. This finishes the proof of (c).\vspace{.1cm}\\(ii) Consider $\chi \in \operatorname{\widetilde{Irr}}(H|\vartheta)$. Let $\lambda$ be an irreducible constituent of $\chi_{\mathcal{A}}.$ We claim that $\lambda \in \operatorname{\widetilde{Lin}}(\mathcal{A}|\vartheta)$. As $A \unlhd H$ and $\vartheta$ is invariant in $H$, by (\cite{IM}, Theorem 6.2), $\vartheta$ is the only irreducible constituent of $\chi_{A}$. However, $\lambda$ being an irreducible constituent of $\chi_{\mathcal{A}}$, it follows that $\lambda_{A}$ is a constituent of $\chi_{A}$. Hence $\vartheta$ is an irreducible constituent of $\lambda_{A}$, and therefore $\operatorname{ker} \vartheta^{\mathcal{A}}\leq \operatorname{ker}\lambda$. But $\operatorname{ker}\vartheta^{\mathcal{A}} = \operatorname{core}_{\mathcal{A}}(\operatorname{ker} \vartheta) = \operatorname{ker} \vartheta$, as $\operatorname{ker} \vartheta$ is normal in $H$. This gives $\operatorname{ker} \vartheta \leq \operatorname{ker} \lambda$. Consequently, $\mathcal{A}/\operatorname{ker} \vartheta$ being abelian, it follows that $\lambda$ is linear and moreover $\lambda_{A} =\vartheta$. We now show that $\operatorname{ker} \lambda^{G}= N$. As $\langle \chi, \lambda^{H}\rangle$, the inner product of $\chi$ with $\lambda^{H}$ is non zero, and $\chi$ is irreducible, we have $\operatorname{ker} \lambda^{H}\leq \operatorname{ker}\chi$. Hence $\operatorname{core}_{G}(\operatorname{ker} \lambda^{H}) \leq \operatorname{core}_{G}(\operatorname{ker}\chi)$, which gives $\operatorname{ker}\lambda^{G} \leq \operatorname{ker}\chi^{G}= N$. Also, $\operatorname{ker}\vartheta \leq \operatorname{ker}\lambda$ implies that $N = \operatorname{core}_{G}(\operatorname{ker} \vartheta)\leq \operatorname{core}_{G}(\operatorname{ker}\lambda)= \operatorname{ker}\lambda^{G}$. Hence the claim follows. Now choose $\varphi \in \operatorname{\mathfrak{Lin}}(\mathcal{A}|\vartheta)$ which lies in the orbit of  $\lambda$. This gives $\sigma \in \operatorname{\mathcal{A}ut}(\mathbb{C}|\vartheta)$ and $x \in H$ such that $\lambda = \sigma\circ \varphi^{x}$. As $\lambda$ is an irreducible constituent of $\chi_{\mathcal{A}},$ it follows that $\varphi^{x}$ is an irreducible constituent of the restriction of $\sigma^{-1}\circ \chi$ to $\mathcal{A}$. However, $\mathcal{A}$ being normal in $H$, from Clifford's theorem (\cite{IM}, Theorem 6.2), it follows that $\varphi$  is an irreducible constituent of the restriction of $\sigma^{-1}\circ \chi$ to $\mathcal{A}$. Consequently, Clifford's correspondence theorem provides $\psi \in \operatorname{{Irr}}(I_{H}(\varphi)|\varphi)$ such that $\sigma^{-1}\circ \chi=\psi^{H}.$ It is easy to check that this $\psi$ belongs to $\operatorname{\widetilde{Irr}}(I_{H}(\varphi)|\varphi)$. Hence (ii) follows. \vspace{.1cm}\\(iii) Suppose $\psi_{1}\in \operatorname{\widetilde{Irr}}(I_{H}(\varphi_{1})|\varphi_{1}),$ $ \psi_{2}\in \operatorname{\widetilde{Irr}}(I_{H}(\varphi_{2})|\varphi_{2})$ are such that \begin{equation}\label{e3} \psi_{2}^{H}= \sigma \circ \psi_{1}^{H}, \end{equation} where $ \sigma \in \operatorname{\mathcal{A}ut}(\mathbb{C})$. By restricting to $\mathcal{A}$, it follows from (\cite{IM}, Theorem 6.2) that \begin{equation}\label{e4} \varphi_{2}=\sigma\circ\varphi_{1}^{h} \end{equation} for some $h \in H$. Therefore, $\varphi_{1}$ and $\varphi_{2}$ lie in the same orbit under the double action and hence $\varphi_{1} = \varphi_{2} = (\varphi~{\rm say})$. In this case, from eqn (\ref{e4}), \begin{equation}\label{e5} \varphi =\sigma\circ\varphi^{h}, \end{equation} which on comparing the kernels yields $ h \in N_{H}(\operatorname{ker}\varphi).$ Now using eqn (\ref{e5}) and the fact that $I_{H}(\varphi) \unlhd N_{H}(\operatorname{ker}\varphi)$, it is easy to see that  $\sigma\circ\psi_{1}^{h} \in \operatorname{\widetilde{Irr}}(I_{H}(\varphi)|\varphi).$ Thus $\sigma\circ\psi_{1}^{h}$ and $\psi_{2}$ both belong to $\operatorname{\widetilde{Irr}}(I_{H}(\varphi)|\varphi)$ and, in view of eqn (\ref{e3}), they are same when induced to $H$. Consequently, the injectivity of the induction map in part (i) implies $\psi_{2} = \sigma\circ\psi_{1}^{h}$. This proves (iii) and completes the proof. $~\Box$

\section{A construction of Shoda pairs} We begin by recalling some basic terminology in graph theory. A graph is a pair $\mathbb{G}=(\mathbb{V}, \mathbb{E})$, where $\mathbb{V}$ is a non-empty set whose elements are termed vertices of $\mathbb{G}$ and $\mathbb{E}$ is a set of unordered pairs of vertices of $\mathbb{G}$. Each element $\{u,v\} \in \mathbb{E}$, where $u,~v \in \mathbb{V}$, is called an edge and is said to join the vertices $u$ and $v$. If $e=\{u,v\} \in \mathbb{E}$, then $e$ is said to be incident with both $u$ and $v$. Further, if $v \in \mathbb{V}$, the degree of $v$ is the number of edges in $\mathbb{E}$ that are incident with $v$. A walk in the graph $\mathbb{G}$ is a sequence of vertices and edges of the form : $v_{1}~\{v_{1},v_{2}\}~v_{2}\cdots \{v_{n-1},v_{n}\}~v_{n}$, where each edge $\{v_{i},v_{i+1}\}$ is incident with the vertices $v_{i}$ and $v_{i+1}$ immediately preceding and succeeding it. A walk is termed path if all the vertices are distinct and is called a cycle if it begins and ends with the same vertex and all other vertices are distinct. A connected graph is the one in which any two vertices are connected by a path. A connected graph which contain no cycles is called a tree.\par A directed graph is a pair $\mathcal{G}=(\mathcal{V}, \mathcal{E})$, where $\mathcal{V}$ is a non-empty set whose elements are termed vertices and $\mathcal{E}$ is a set of ordered pairs of vertices of $\mathcal{V}$. In a directed graph, an edge $e=(u,v)$ is said to be incident out of $u$ and incident into $v$. The terminology of directed walk and directed path is same as that in graph but now the edges are directed in the same direction. In a directed graph, the number of edges incident out of $v$ is called out-degree of $v$ and the number of edges incident into $v$ is called in-degree of $v$, where $v \in \mathcal{V}$. The vertices of in-degree $0$ are called source and those of out-degree $0$ are called sink. In a directed graph, there is an obvious underlying undirected graph whose vertex set is same as that of the directed graph and there is an edge $\{u,v\}$ if either $(u,v)$ or $(v,u)$ is an edge in the directed graph. A directed graph is called a directed tree if its underlying undirected graph is a tree. The sink vertices of a directed tree are termed as its leaves. A directed tree is called a rooted directed tree if it has a unique source. The unique source of a rooted directed tree is called its root. \par We now proceed with the construction of Shoda pairs. Let $G \in \mathcal{C}$ and let $N$ be a normal subgroup of $G$. Consider the directed graph $\mathcal{G}=(\mathcal{V}, \mathcal{E})$ whose vertex set $\mathcal{V}$ consist of all $N$-linear character triples of $G$ and there is an edge $(u,v) \in \mathcal{E}$ if $v$ is a direct Clifford correspondent(d.c.c.) of $u$. Clearly $(G,N,1_{N}) \in \mathcal{V}$, where $1_{N}$ is the character of $N$ which takes constant value $1$. Let $\mathcal{V}_{N}$ be the set of those vertices $v \in \mathcal{V}$ for which there is a directed path from $(G,N,1_{N})$ to $v$. Let $\mathcal{E}_{N}$ be the set of ordered pairs $(u,v)\in \mathcal{E}$ with $u,v \in \mathcal{V}_{N}$. Then $\mathcal{G}_{N}=(\mathcal{V}_{N},\mathcal{E}_{N})$ is a directed subgraph of $\mathcal{G}$. Observe that any vertex $(H, A, \vartheta)$ of $\mathcal{G}_{N}$ with $H=A$ is a sink vertex. \begin{theorem}\label{t2} Let $G \in \mathcal{C}$ and $\mathcal{N}$ the set of all normal subgroups of $G$. \begin{description}  \item[(i)] For $N \in \mathcal{N}$, the following hold: \begin{description}  \item[(a)]  $\mathcal{G}_{N}$ is a rooted directed tree with $(G,N,1_{N})$ as its root; \item[(b)] the leaves of $\mathcal{G}_{N}$ of the type $(H, H, \vartheta)$ correspond to Shoda pairs of $G$. More precisely, if $(H, H, \vartheta)$ is a leaf of $\mathcal{G}_{N}$, then  $(H,\operatorname{ker}\vartheta)$ is a Shoda pair of $G$.\end{description} \item[(ii)] If  $({H'}, {K'})$ is any Shoda pair of $G$, then  there is a leaf $(H,H,\vartheta)$ of $\mathcal{G}_{N}$, where $N=\operatorname{core}_{G}(K')$,  such that $({H'}, {K'})$ and $(H,\operatorname{ker}\vartheta)$ realize the same primitive central idempotent of $\mathbb{Q}G$. \end{description} \end{theorem}\vspace{.2cm} \par To prove the theorem, we need some preparation.
\begin{lemma}\label{l1}  For each vertex $v$ of $\mathcal{G}_{N}$, there is a unique directed path from $(G,N,1_{N})$ to $v$.  \end{lemma}
 {\bf Proof.}
 Let $v_{1}~(v_{1},v_{2})~v_{2}\cdots (v_{n-1}, v_{n})~v_{n}$ and $v_{1}'~(v_{1}',v_{2}')~v_{2}'\cdots (v_{m-1}', v_{m}')~v_{m}'$ be two directed paths from $v_{1}= v_{1}'=(G,N,1_{N})$ to $v_{n}=v_{m}'=v=(H, A, \vartheta)$. Assume that $m \leq n$. We claim that $ v_{i}=v_{i}'~{\rm for}~  1 \leq i \leq m$. We'll prove it by induction on $i$. For $i=1$, we already have $v_{1}=v_{1}' = (G,N,1_{N}) $. Assume that $v_{i}=v_{i}'$ for some $i<m$. Write $v_{j}=(H_{j}, A_{j}, \vartheta_{j })$ and $v_{j}'=(H_{j}', A_{j}', \vartheta_{j}')$ for $1 \leq j \leq m$. From the construction of d.c.c., we have $A_{i+1} = \mathcal{A}_{(H_{i}, A_{i}, \vartheta_{i})}$ and $A_{i+1}' = \mathcal{A}_{(H_{i}', A_{i}', \vartheta_{i}')}$. As $(H_{i}, A_{i}, \vartheta_{i})= (H_{i}', A_{i}', \vartheta_{i}')$,  it follows immediately that $A_{i+1}=A_{i+1}'$. Now $\vartheta_{i+1}$ (resp. $\vartheta_{i+1}'$) being the restriction of $\vartheta$ to $A_{i+1}$ (resp. $A_{i+1}'$) yields that $\vartheta_{i+1}= \vartheta_{i+1}'$. Further, as $H_{i+1}= I_{H_{i}}(\vartheta_{i+1})$ and $H_{i+1}'= I_{H_{i}'}(\vartheta_{i+1}')$, it follows that $H_{i+1}= H_{i+1}'$. This proves the claim, which as a consequence implies that both $v_{m}$ and $v_{n}$ are equal to $v$. This is not possible if $m < n$ as no two vertices in a path are same. $~\Box$

 \begin{lemma}\label{l2} The following statements hold for $\mathcal{G}_{N}$: \begin{description}\item[(i)] If $(u_{1}, v_{1})$ and $(u_{2}, v_{2})$ are edges of $\mathcal{G}_{N}$ with $v_{1}=v_{2}$, then $u_{1}=u_{2}$; \item[(ii)] If $v_{1}~\{v_{1},v_{2}\}~v_{2}\cdots \{v_{n-1},v_{n}\}~v_{n}$ is a path in the underlying undirected graph of $\mathcal{G}_{N}$, then there is a unique $j$, $1\leq j\leq n$, with the following: \begin{description} \item[(a)] $(v_{i+1},v_{i}) \in \mathcal{E}_{N}$  for $1\leq i < j$;  \item[(b)] $(v_{i},v_{i+1})\in \mathcal{E}_{N}$ for $ j \leq i < n$. \end{description} \end{description}\end{lemma}
 {\bf Proof.} (i) is a consequence of Lemma \ref{l1} and (ii) follows immediately from (i). It may be mentioned that in the lemma, (a) is empty if $j=1$, and (b) is empty if $j=n$. $~\Box$
\begin{lemma}\label{l3} The underlying undirected graph of $\mathcal{G}_{N}$ is a tree.\end{lemma}
{\bf Proof.} Let $\mathbb{G}_{N}=(\mathbb{V}_{N},\mathbb{E}_{N})$ be the underlying undirected graph of $\mathcal{G}_{N}=(\mathcal{V}_{N},\mathcal{E}_{N}).$ Clearly $\mathbb{G}_{N}$ is a connected graph. To show that $\mathbb{G}_{N}$ is a tree, it is enough to prove that $\mathbb{G}_{N}$ is disconnected after removing an edge (see \cite{ND}, Theorem 3-5). Let $e=\{u,v\} \in \mathbb{E}_{N}$ and let $\mathbb{E}'_{N} = \mathbb{E}_{N} \setminus \{e\}$. We need to show that  $\mathbb{G}'_{N}=(\mathbb{V}_{N},\mathbb{E}'_{N})$ is disconnected. If not, then there is a path \begin{equation} \label{e6} v_{1}~\{v_{1},v_{2}\}~v_{2}\cdots \{v_{n-1},v_{n}\}~v_{n}\end{equation} in $\mathbb{G}_{N}'$, where $v_{1}=u$ and $v_{n}=v$. As $\{u,v\} \in \mathbb{E}_{N}$, either $(u,v)$ or $(v,u)$ belongs to $\mathcal{E}_{N}$. Suppose $(u, v) \in \mathcal{E}_{N}$. As the path given in eqn (\ref{e6}) is also a path in $\mathbb{G}_{N}$, there is a $j$, $1 \leq j \leq n$ so that part (ii) of Lemma \ref{l2} holds. If $j<n$, then $v_{n}=v$ is a d.c.c. of $v_{n-1}$, which by Lemma \ref{l2}(i) implies $v_{n-1}=u$. Consequently $\{u,v\}= \{v_{n-1}, v_{n}\} \in \mathbb{E}'_{N},$ which is not so. If $j=n$, then $v_{i}$ is a d.c.c. of $v_{i+1}$ for all $1\leq i < n$. Hence if $v_{i}=(H_{i}, A_{i}, \vartheta_{i})$, then from Theorem \ref{t1}, $A_{n} \lneq A_{n-1} \lneq \cdots \lneq A_{1}$. However $(u,v) =(v_{1},v_{n}) \in \mathcal{E}_{N}$ implies $A_{1}\lneq A_{n}$, a contradiction. Using similar arguments, it follows that $(v,u) \in \mathcal{E}_{N}$ is not possible. This proves the lemma.
$~\Box$ \vspace{.2cm}\par We are now ready to prove the theorem. Recall from (\cite{TY}, Proposition 1.1) that if $\chi \in \operatorname{Irr}G$, then $e_{\mathbb{Q}}(\chi) = \frac{\chi(1)}{|G|}\sum_{\sigma \in \operatorname{Gal}(\mathbb{Q}(\chi)/\mathbb{Q})} \sum_{g \in G}\sigma(\chi(g))g^{-1}$ is a primitive central idempotent of $\mathbb{Q}G$, where $\operatorname{Gal}(\mathbb{Q}(\chi)/\mathbb{Q})$ is the Galois group of $\mathbb{Q}(\chi)$ over $\mathbb{Q}$. \vspace{.2cm} \\ \noindent{\bf Proof of Theorem \ref{t2}} (i) By Lemma \ref{l3}, $\mathcal{G}_{N}$ is a directed tree. If the in-degree of $(G,N,1_{N})$ is non zero then there is a $N$-linear character triple $(H,A,\vartheta)$ of $G$ such that $(G,N,1_{N})$ is a d.c.c. of $(H,A,\vartheta)$. Hence by Theorem \ref{t1}(i)(a) $A\lneq N$. However $(H,A,\vartheta)$ being $N$-linear character triple of $G$, we have $\operatorname{ker}\vartheta^{G}=N$, which gives $N\leq A$, a contradiction. This proves that $(G,N,1_{N})$ is a source of $\mathcal{G}_{N}$. Next if $(H,A,\vartheta)$ is a vertex of $\mathcal{G}_{N}$ different from $(G,N,1_{N})$, then there is a directed path from $(G,N,1_{N})$ to $(H,A,\vartheta)$, which implies that the in-degree of $(H,A,\vartheta)$ is non zero. Hence $(H,A,\vartheta)$ can't be a source. This proves (a). To prove (b), consider a leaf of $\mathcal{G}_{N}$ of the type $(H,H,\vartheta)$. Let $v_{1}~(v_{1},v_{2})~v_{2}\cdots (v_{n-1},v_{n})~v_{n}$ be the directed path from $v_{1} = (G,N,1_{N})$ to $v_{n}=(H,H,\vartheta)$. Let $v_{i}=(H_{i},A_{i},\vartheta_{i})$, $1\leq i \leq n$. As $(H_{i+1},A_{i+1},\vartheta_{i+1})$ is a d.c.c. of $(H_{i},A_{i},\vartheta_{i})$, from Theorem \ref{t1}
 \begin{equation} \label{e7} \psi^{H_{i}} \in \operatorname{\widetilde{Irr}}(H_{i}|\vartheta_{i})~{\rm for~all}~ \psi \in \operatorname{\widetilde{Irr}}(H_{i+1}|\vartheta_{i+1}) .\end{equation} This is true for all $1\leq i < n$. Observe that $\operatorname{\widetilde{Irr}}(H_{n}|\vartheta_{n}) = \{\vartheta_{n} \} = \{\vartheta \}$, as $H_{n}=A_{n}=H$. The repeated application of eqn (\ref{e7}) when $\psi$ is $\vartheta$ implies that \begin{equation}\label{e8} \vartheta^{H_{i}} \in \operatorname{\widetilde{Irr}}(H_{i}| \vartheta_{i})~{\rm  for}~ 1\leq i \leq n. \end{equation} For $i=1$, the above equation, in particular, gives  $\vartheta^{G} \in \operatorname{Irr} G$. This proves (b).\\ (ii) Let $(H',K')$ be a Shoda pair of $G$ and let $N = \operatorname{core}_{G}(K')$. Let $\psi$ be a linear character on $H'$ with kernel $K'$. We claim that there is a leaf $(H,H,\vartheta)$ of $\mathcal{G}_{N}$ such that $(H',K')$ and $(H,\operatorname{ker}\vartheta)$ realize the same primitive central idempotent of $\mathbb{Q}G.$ From the definition of the primitive central idempotent of $\mathbb{Q}G$ realized by a Shoda pair given in \cite{OdRS04}, it follows that the primitive central idempotent of $\mathbb{Q}G$ realized by $(H', K')$ is $e_{\mathbb{Q}}(\psi^{G})$ and that realized by $(H,\operatorname{ker}\vartheta)$ is $e_{\mathbb{Q}}(\vartheta^{G})$. Let $\chi =\psi^{G}$. By Shoda's theorem, $\chi \in \operatorname{Irr} G.$ Also, $\operatorname{ker}\chi = N$. If $N=G$, then $(H',K')$ is clearly $(G,G)$. Also in this case, $\mathcal{G}_{N}$ is just the vertex $(G,G,1_{G})$ which corresponds to the Shoda pair $(G,G)=(H',K').$ Thus we may assume that $N\neq G$. Denote $\chi$ by $\chi_{1}$. As $\chi_{1}\in \operatorname{\widetilde{Irr}}(G|1_{N})$, by Theorem \ref{t1}, there is a d.c.c., $(I_{G}(\varphi),\mathcal{A}_{(G,N,1_{N})},\varphi)$ of $(G,N,1_{N})$, $\chi_{2} \in \operatorname{\widetilde{Irr}}(I_{G}(\varphi)|\varphi)$ and $\sigma_{1} \in \operatorname{\mathcal{A}ut}(\mathbb{C})$ such that $\chi_{1}=\sigma_{1}\circ\chi_{2}^{G}.$ Put $(H_{1},A_{1},\vartheta_{1})=(G,N,1_{N})$, $(H_{2},A_{2},\vartheta_{2})= (I_{G}(\varphi),\mathcal{A}_{(G,N,1_{N})},\varphi)$. If $H_{2}= A_{2}$, stop. If not, again by Theorem \ref{t1}, there is a d.c.c. $(H_{3},A_{3},\vartheta_{3})$ of $(H_{2},A_{2},\vartheta_{2})$, $\chi_{3} \in \operatorname{\widetilde{Irr}}(H_{3}|\vartheta_{3})$, $\sigma_{2} \in \operatorname{\mathcal{A}ut}(\mathbb{C})$ such that \begin{equation} \label{e9}  \chi_{2}=\sigma_{2}\circ\chi_{3}^{H_{2}}.\end{equation} Moreover, if this case arises, then by Theorem \ref{t1}(i)(a), $N= A_{1}\lneq A_{2}\lneq A_{3} \leq G.$ Again if $H_{3}= A_{3}$ stop, otherwise continue. This process of continuing must stop after finite number of steps as at nth step, there is an ascending chain $$N= A_{1}\lneq A_{2}\lneq \cdots \lneq A_{n} \leq G.$$ Suppose the process stops at nth step. Then we have character triples $(H_{i},A_{i},\vartheta_{i})$, $1\leq i\leq n$, with  $H_{n}=A_{n}$, $\chi_{i} \in \operatorname{\widetilde{Irr}}(H_{i}|\vartheta_{i})$ and $\sigma_{i} \in \operatorname{\mathcal{A}ut}(\mathbb{C})$ such that \begin{equation} \label{e10} \chi_{i}=\sigma_{i}\circ\chi_{i+1}^{H_{i}},~~~1\leq i < n. \end{equation} The above equation yields $$\chi =\sigma\circ\chi_{n}^{G}, {\rm~~~where~} \sigma =\sigma_{1}\circ\sigma_{2}\circ \cdots \circ\sigma_{n}.$$ As $H_{n}=A_{n}$, we have $\operatorname{\widetilde{Irr}}(H_{n}|\vartheta_{n})=\{\vartheta_{n}\}$, and hence $\chi_{n}=\vartheta_{n}$. This gives $\chi =\sigma\circ\theta_{n}^{G}$. Consequently $e_{\mathbb{Q}}(\chi)=e_{\mathbb{Q}}(\vartheta_{n}^{G})$ and hence $(H',K')$ and $(H_{n},\operatorname{ker}\vartheta_{n}$) realize the same primitive central idempotent of $\mathbb{Q}G$. This proves the claim and completes the proof of theorem. $~\Box$ \vspace{.3cm} \par For $N \in \mathcal{N}$, denote by $\mathcal{L}_{N}$ the set of leaves of $\mathcal{G}_{N}$ of type $(H,H,\vartheta)$. Let $\mathcal{S}_{N}$ be the set of Shoda pairs of $G$ corresponding to the leaves in $\mathcal{L}_{N}$. We have shown in Theorem \ref{t2} that if $G \in \mathcal{C}$, then the mapping from $\bigcup_{N \in \mathcal{N}}S_{N} $ to the set of primitive central idempotents of $\mathbb{Q}G$ realized by the Shoda pairs of $G$ is surjective. In other words, $\bigcup_{N \in \mathcal{N}}\mathcal{S}_{N}$ is a complete set of Shoda pairs of $G$. We now begin to investigate whether this set of Shoda pairs of $G$ is irredundant, i.e., this map is injective. \par Two Shoda pairs $(H_{1},K_{1})$ and $(H_{2},K_{2})$ of $G$ are said to be equivalent if they realize the same primitive central idempotent of $\mathbb{Q}G.$ \begin{lemma}\label{l4} If $N$ and $N'$ are distinct normal subgroups of $G$, then the Shoda pair corresponding to a leaf in $\mathcal{L}_{N}$ can not be equivalent to that in $\mathcal{L}_{N'}$.\end{lemma}{\bf Proof.} Let $(H,H,\vartheta) \in \mathcal{L}_{N}$ and $(H',H',\vartheta') \in \mathcal{L}_{N'}$. Suppose $(H,\operatorname{ker}\vartheta)$ is equivalent to $(H',\operatorname{ker}\vartheta')$. Then $e_{\mathbb{Q}}(\vartheta^{G})=e_{\mathbb{Q}}(\vartheta'^{G})$, which gives $\vartheta'^{G}=\sigma\circ\vartheta^{G}$ for some $\sigma \in \operatorname{\mathcal{A}ut}(\mathbb{C}).$ Consequently, $\operatorname{ker}\vartheta'^{G}=\operatorname{ker}\sigma\circ\vartheta^{G}=\operatorname{ker}\vartheta^{G}$. Hence, $N=N'$. $~\Box$\vspace{.3cm} \par We next examine if distinct leaves in $\mathcal{L}_{N}$, for a fixed normal subgroup $N$ of $G$, can correspond to equivalent Shoda pairs. For this purpose, we need to fix some terminology. If $(H,H,\vartheta) \in \mathcal{L}_{N}$ and $v_{1}~(v_{1},v_{2})~v_{2}\cdots (v_{n-1},v_{n})~v_{n}$ is the directed path from $v_{1}=(G,N,1_{N})$ to $v_{n}=(H,H,\vartheta)$, then we call $n$ to be the height of $(H,H,\vartheta)$ and term $v_{i}$ as the $i^{th}$ node of $(H,H,\vartheta)$, $1\leq i\leq n.$ It may be noted that if $v_{i}=(H_{i},A_{i},\vartheta_{i})$, then from eqn (\ref{e8}), it follows that $\vartheta^{H_{i}} \in \operatorname{\widetilde{Irr}}(H_{i}|\vartheta_{i})$ for all $1 \leq i\leq n.$\vspace{.2cm}\\${\rm \mathbf{Definition}}$ Let $(H,H,\vartheta) \in \mathcal{L}_{N}$ be of height $n$ with $(H_{i},A_{i},\vartheta_{i})$ as its $i^{th}$ node, $1\leq i\leq n$. We call $(H,H,\vartheta)$ to be good if the following holds for all $1< i\leq n$: given $x \in N_{H_{i-1}}(\operatorname{ker}\vartheta_{i})$, there exist $\sigma \in \operatorname{\mathcal{A}ut}(\mathbb{C})$ such that $(\vartheta^{H_{i}})^{x}=\sigma\circ\vartheta^{H_{i}}.$\vspace{.2cm}\begin{remark}\label{r1} It may be noted that, for any normal subgroup $N$ of $G$, the leaves of $\mathcal{L}_{N}$ of height $2$ are always good.\end{remark}\begin{lemma}\label{l5} If $N \in \mathcal{N}$ and $(H,H,\vartheta) \in \mathcal{L}_{N} $ is good, then its corresponding Shoda pair can not be equivalent to that of other leaves in $\mathcal{L}_{N}.$ \end{lemma}{\bf Proof.} Let $(H,H,\vartheta)$ and $(H',H',\vartheta') \in \mathcal{L}_{N}$ be distinct and let $(H,H,\vartheta)$ be good. Let the height of $(H,H,\vartheta)$ and $(H',H',\vartheta')$ be $n$ and $n'$ respectively. Let $(H_{i},A_{i},\vartheta_{i})$ be the $i^{th}$ node of $(H,H,\vartheta)$, $1\leq i\leq n$, and $(H_{j}',A_{j}',\vartheta_{j}')$ the $j^{th}$ node of $(H',H',\vartheta')$, $1\leq j\leq n'$. We have\begin{equation}\label{e11} \vartheta^{H_{i}} \in \widetilde{\operatorname{Irr}}(H_{i}|\vartheta_{i}),~1\leq i\leq n\end{equation} and  \begin{equation}\label{e12} \vartheta'^{H_{j}'} \in \widetilde{\operatorname{Irr}}(H_{j}'|\vartheta_{j}'),~1\leq j\leq n'.\end{equation} Let $k$ be the least positive integer such that $(H_{k},A_{k},\vartheta_{k})\neq (H_{k}',A_{k}',\vartheta_{k}').$ Clearly, $k \neq 1$. By the definition of $k$, $(H_{k},A_{k},\vartheta_{k})$ and $(H_{k}',A_{k}',\vartheta_{k}')$ are distinct d.c.c.'s of $(H_{k-1},A_{k-1},\vartheta_{k-1})$. The construction of d.c.c. yields $A_{k}=A_{k'}$. In view of eqns (\ref{e11}), (\ref{e12}), it follows from Theorem \ref{t1}(iii) that \begin{equation}\label{e13}\vartheta'^{H_{k-1}}\neq \sigma\circ\vartheta^{H_{k-1}}\end{equation} for any $\sigma \in \operatorname{\mathcal{A}ut}(\mathbb{C})$. We now show  a contradiction to eqn (\ref{e13}), if $(H,\operatorname{ker}\vartheta)$ and $(H',\operatorname{ker}\vartheta')$ are equivalent. We have that $H_{i}=H_{i}'$ and $\vartheta_{i}=\vartheta_{i}'$ for $1\leq i\leq k-1$. Also $H_{1}=H_{1}'=G$. From eqns (\ref{e11}) and (\ref{e12}), $\vartheta^{H_{2}}$ and $\vartheta'^{H_{2}}$ belong to $\widetilde{\operatorname{Irr}}(H_{2}|\vartheta_{2}).$ If $(H,\operatorname{ker}\vartheta)$ and $(H',\operatorname{ker}\vartheta')$ are equivalent, then $e_{\mathbb{Q}}(\vartheta'^{G})= e_{\mathbb{Q}}(\vartheta^{G})$, which gives $$\vartheta'^{H_{1}}=\sigma_{1}\circ\vartheta^{H_{1}}, {\rm~~~ for~some ~}\sigma_{1} \in \operatorname{\mathcal{A}ut}(\mathbb{C}).$$ Hence, by Theorem \ref{t1}(iii), there exist $x \in N_{H_{1}}(\operatorname{ker}\vartheta_{2})$ such that \begin{equation}\label{e14}\vartheta'^{H_{2}}=\sigma_{1}\circ(\vartheta^{H_{2}})^{x}.\end{equation} Also, $(H,H,\vartheta)$ being good, there exists $\tau \in \operatorname{\mathcal{A}ut}(\mathbb{C})$ such that \begin{equation}\label{e15}(\vartheta^{H_{2}})^{x}=\tau\circ\vartheta^{H_{2}}.\end{equation} Let $\sigma_{2}=\sigma_{1}\circ\tau$. From eqns (\ref{e14}) and (\ref{e15}), $$\vartheta'^{H_{2}}=\sigma_{2}\circ\vartheta^{H_{2}}.$$ Now repeating this process with $H_{1}$, $H_{2}$ replaced by $H_{2}$, $H_{3}$ respectively, we obtain that $$\vartheta'^{H_{3}}=\sigma_{3}\circ\vartheta^{H_{3}},~{\rm for~some}~\sigma_{3} \in \operatorname{\mathcal{A}ut}(\mathbb{C}).$$ Continuing this process, we obtain after $k-1$ steps that $$\vartheta'^{H_{k-1}}=\sigma_{k-1}\circ\vartheta^{H_{k-1}},~{\rm for~some}~\sigma_{k-1} \in \operatorname{\mathcal{A}ut}(\mathbb{C}).$$ This contradicts eqn (\ref{e13}) and completes the proof. $~\Box$
\begin{theorem}\label{t3} For $G \in \mathcal{C}$, $\bigcup_{N \in \mathcal{N}}\mathcal{S}_{N}$ is a complete irredundant set of Shoda pairs of $G$ if, and only if, the leaves in $\mathcal{L}_{N}$ are good for all $N \in \mathcal{N}.$\end{theorem} {\bf Proof.} Suppose $\bigcup_{N \in \mathcal{N}}\mathcal{S}_{N}$ is a complete irredundant set of Shoda pairs of $G \in \mathcal{C}$. Let $N \in \mathcal{N}$. Let $(H,H,\vartheta)\in \mathcal{L}_{N}$ be of height $n$ with $(H_{i},A_{i},\vartheta_{i})$ its $i^{th}$ node, $ 1 \leq i \leq n$. We'll show that $(H,H,\vartheta)$ is good. Let $2\leq i\leq n$ and $x \in N_{H_{i-1}}(\operatorname{ker}\vartheta_{i}).$ By Theorem \ref{t1}(i)(c), $(\vartheta^{H_{i}})^{x}\in \widetilde{\operatorname{Irr}}(H_{i}|\vartheta_{i}).$  Proceeding as in the proof of part (ii) of Theorem \ref{t2}, there exists $(H',H',\vartheta')\in \mathcal{L}_{N}$ such that\begin{equation}\label{e16}
 (\vartheta^{H_{i}})^{x}=\sigma\circ\vartheta'^{H_{i}}, {\rm~for~some~}\sigma \in \operatorname{\mathcal{A}ut}(\mathbb{C}).\end{equation} Inducing to $G$, we get $\vartheta^{G}=\sigma\circ\vartheta'^{G}$, which gives $\vartheta=\vartheta'$ and $H=H'$, as $\bigcup_{N \in \mathcal{N}}\mathcal{S}_{N}$ is complete and irredundant. Consequently, eqn (\ref{e16}) implies that $(H,H,\vartheta)$ is good. This proves `only if' statement. The `if' statement follows from Theorem \ref{t2} and Lemma \ref{l5}.  $~\Box$  \vspace{.1cm} \par As a consequence, the above theorem yields the following result proved in (\cite{BM}, Theorem 1(i)):
  \begin{cor}\label{c0} If $G$ is a normally monomial group, then $\bigcup_{N \in \mathcal{N}}\mathcal{S}_{N}$ is a complete irredundant set of Shoda pairs of $G$.\end{cor}\noindent \textbf{Proof.} It is enough to show that all the character triples in $\bigcup_{N \in \mathcal{N}}\mathcal{L}_{N}$ are good. Let $(H,H,\vartheta) \in \mathcal{L}_{N}$, where $N \in \mathcal{N}$. If $N=G$, then $(H,H,\vartheta)=(G,G,1_{G})$, which is clearly good. Assume $N\neq G$. As $\vartheta^{G} \in \widetilde{\operatorname{Irr}}(G|1_{N})$, $\operatorname{ker}\vartheta^{G}=N$. Hence, by (\cite{HowIII}, Lemma 2.2), $\vartheta^{G}(1)=[G:\mathcal{A}_{(G,N,1_{N})}]$. However, $\vartheta^{G}(1)=[G:H]$. Therefore, we have \begin{equation}\label{e0} [G:\mathcal{A}_{(G,N,1_{N})}]=[G:H]\end{equation} Let the height of $(H,H,\vartheta)$ be $n$ and let $(H_{i},A_{i},\vartheta_{i})$ be its $i^{th}$ node. As $(H_{2},A_{2},\vartheta_{2})$ is d.c.c. of $(H_{1},A_{1},\vartheta_{1})=(G,N,1_{N})$, we have $A_{2}=\mathcal{A}_{(G,N,1_{N})}$. Also in view of Theorem \ref{t1}(i)(a), we have $N=A_{1}\lneq A_{2}\lneq \cdots \lneq A_{n}=H$. This gives $\mathcal{A}_{(G,N,1_{N})}\lneq H$, if $n>2$, in which case eqn (\ref{e0}) can't hold. Hence we must have $n=2$, which in view of remark \ref{r1}, implies that $(H,H,\vartheta)$ is good. $~\Box$
  \begin{remark}\label{r0} Later from remark \ref{r2}, it will follow that if $G$ is a normally monomial group, then the Shoda pairs in $\bigcup_{N\in \mathcal{N}}\mathcal{S}_{N}$ are strong Shoda pairs of $G$.\end{remark}

 \section{Idempotents from Shoda pairs} We continue to use the notation developed in the previous section. Given a group $G \in \mathcal{C}$, we have shown in the previous section that any Shoda pair of $G$ is equivalent to $(H, \operatorname{ker} \vartheta)$, where $(H, H, \vartheta)$ is a character triple in $\bigcup_{N \in \mathcal{N}}\mathcal{L}_{N}$, and it realizes the primitive central idempotent $e_{\mathbb{Q}}(\vartheta^{G})$ of $ \mathbb{Q}G$. In this section, we examine the expression of $e_{\mathbb{Q}}(\vartheta^{G})$ in terms of $e(G, H, K)$, where $K=\operatorname{ker} \vartheta$. In \cite{OdRS04}, Olivieri, del R{\'{\i}}o and Sim{\'o}n proved that $e_{\mathbb{Q}}(\vartheta^{G})=\alpha_{(G,H,K)} e(G,H,K)$, where $\alpha_{(G,H,K)}=\frac{[\operatorname{Cen}_{G}(\varepsilon(H,K)):H]}{[\mathbb{Q}(\vartheta):\mathbb{Q}(\vartheta^{G})]}$. The following theorem gives a new character free expression of $\alpha_{(G,H,K)}$ and also provides a necessary and sufficient condition for $\alpha_{(G,H,K)}= 1$. It may be mentioned that $\alpha_{(G,H,K)}=1$ is a necessary condition for $(H,K)$ to be a strong Shoda pair of $G$. \begin{theorem}\label{t4} Let $G \in \mathcal{C}$ and $N$ a normal subgroup of $G$. Let $(H,H,\vartheta) \in \mathcal{L}_{N}$ be of height $n$ with $(H_{i},A_{i},\vartheta_{i})$  as its $i^{th}$ node, $1\leq i\leq n.$ Let $K= \operatorname{ker}\vartheta.$ Then \begin{description}\item[(i)] $\alpha_{(G,H,K)}=\frac{[\operatorname{Cen}_{G}(\varepsilon(H,K)):\operatorname{Cen}_{H_{n-1}}(\varepsilon(H,K))]}{\prod_{2\leq i\leq n-1}[\operatorname{Cen}_{H_{i-1}}(e(H_{i},H,K)):H_{i}]}$; \item[(ii)]$\alpha_{(G,H,K)}=1$ if, and only if,~$\operatorname{Cen}_{H_{i-1}}(e(H_{i},H,K))= \operatorname{Cen}_{H_{i-1}}(\varepsilon(H,K))H_{i}$ for\linebreak all $2\leq i\leq n-1$;\item[(iii)]if $(H,H,\vartheta)$ is good, then in the above statements $\operatorname{Cen}_{H_{i-1}}(e(H_{i},H,K))$ can \linebreak be replaced by $N_{H_{i-1}}(\operatorname{ker}\vartheta_{i})$ for all $2\leq i\leq n-1$. \end{description}\end{theorem}\vspace{.2cm} \par We first prove the following:
\begin{lemma}\label{l6} Let $G$ be a finite group and  $K\unlhd H\leq G $ with $H/K$ cyclic. Let $A\unlhd H$ and $D=K\cap A.$ Then $$\varepsilon(A,D)=\varepsilon(H,K)+e$$ for some central idempotent $e$ of $\mathbb{Q}H$ orthogonal to $\varepsilon(H,K).$ \end{lemma}{\bf Proof.} Let $\psi$ be a linear character on $H$ with kernel $K$ and let $\varphi= \psi_{A}$. By (\cite{IM}, Corollary 6.17), we have $$\varphi^{H}=\displaystyle\sum_{\chi \in \operatorname{Irr}H/A} \langle \varphi^{H},  \chi\psi \rangle \chi\psi.$$ As $\varphi$ is invariant in $H$, we have $(\chi\psi)_{A} = \langle \varphi,  (\chi\psi)_{A} \rangle \varphi= \langle \varphi^{H},  \chi\psi \rangle \varphi$, and hence $\chi\psi(1)= \langle \varphi^{H},  \chi\psi \rangle \varphi(1) = \langle \varphi^{H},  \chi\psi \rangle$. Therefore, $$\varphi^{H}=\displaystyle\sum_{\chi \in \operatorname{Irr}H/A} (\chi\psi)(1) \chi\psi.$$ Let $n$=[$A$:$D$] and $\sum$ the collection of all the irreducible constituents of $(\varphi^{i})^{H},\linebreak 1\leq i\leq n$ with $(i,n)=1$. In other words, $$\sum = \{\chi\psi^{i}~|~\chi \in \operatorname{Irr}H/A,~1\leq i\leq n,~(i,n)=1\}.$$ Consider the natural action of $\operatorname{\mathcal{A}ut}(\mathbb{C})$ on $\sum$. Under this action, let $\chi_{1},~\chi_{2},\cdots,~\chi_{r}$ be the representatives of distinct orbits with $\chi_{1}=\psi$. It can be checked that $\operatorname{orb}(\chi_{i})$, the orbit of $\chi_{i}$, is given by $\{\sigma\circ\chi_{i}~|~\sigma \in \operatorname{Gal}(\mathbb{Q}(\chi_{i})/\mathbb{Q})\}$. Also, we have  $$\displaystyle\sum_{1\leq i\leq n,~(i,n)=1}(\varphi^{i})^{H}= \displaystyle\sum_{1\leq i\leq r}\sum_{\sigma \in \operatorname{Gal}(\mathbb{Q}(\chi_{i})/\mathbb{Q})} \chi_{i}(1) \sigma\circ\chi_{i}.$$ Consequently, $$\begin{array}{lll} \displaystyle\sum_{1\leq i\leq r}e_{\mathbb{Q}}(\chi_{i})&=&\displaystyle\sum_{1\leq i\leq r}\frac{1}{|H|}\sum_{h \in H}\sum_{\sigma \in \operatorname{Gal}(\mathbb{Q}(\chi_{i})/\mathbb{Q})}\chi_{i}(1)\sigma(\chi_{i}(h))h^{-1}\\ &=&\displaystyle\frac{1}{|H|}\sum_{h \in H} \displaystyle\sum_{1\leq i\leq n,~(i,n)=1}(\varphi^{i})^{H}(h)h^{-1}\\&=&\displaystyle\frac{1}{|H|}\sum_{h \in H}\sum_{\tau \in \operatorname{Gal}(\mathbb{Q}(\varphi)/\mathbb{Q})}(\tau\circ\varphi)^{H}(h)h^{-1} \\&=&\displaystyle\frac{1}{|A|}\sum_{a \in A}\sum_{\tau \in \operatorname{Gal}(\mathbb{Q}(\varphi)/\mathbb{Q})}(\tau\circ\varphi)(a)a^{-1}\\&=&\varepsilon(A,D).\end{array}$$ As $e_{\mathbb{Q}}(\chi_{1})= e_{\mathbb{Q}}(\psi)=\varepsilon(H,K)$, the result follows. $~\Box$\vspace{.2cm} \par The following proposition is crucial in the proof of Theorem \ref{t4}. \begin{prop}\label{p1} Let $G$ be a finite group. Let $K\unlhd H\leqslant G$ with $H/K$ cyclic and $\psi$ a linear character on $H$ with kernel $K$. Suppose that there is a normal subgroup $A$ of $G$ and a subgroup $L$ of $G$ such that $A \leq H \leq L$, $\psi_{A}$ is invariant in $L$, and $L\unlhd N_{G}(\operatorname{ker}\psi_{A})$. Then the following hold:\begin{description}\item[(i)] if $\psi^{L} \in \operatorname{Irr}L$, then the distinct $G$-conjugates of $e_{\mathbb{Q}}(\psi^{L})$ are mutually orthogonal;\item[(ii)]if $\psi^{G} \in \operatorname{Irr}G$, then $e_{\mathbb{Q}}(\psi^{G})$ is the sum of all distinct $G$-conjugates of $e_{\mathbb{Q}}(\psi^{L}).$ Furthermore, $$\alpha_{(G,H,K)}=\alpha_{(L,H,K)}\frac{[\operatorname{Cen}_{G}(\varepsilon(H,K)):\operatorname{Cen}_{L}(\varepsilon(H,K))]}
{[\operatorname{Cen}_{G}(e(L,H,K)):L]}.$$\end{description} \end{prop}\noindent{\bf Proof.} (i) Denote $\operatorname{ker}\psi_{A}$ by $D$. First of all, we will
show that \begin{equation}\label{e17}e_{\mathbb{Q}}(\psi^{L})e_{\mathbb{Q}}(\psi^{L})^{g}=0,~{\rm~if~}g\not\in N_{G}(D).\end{equation}For this, it is enough to prove that \begin{equation}\label{e18}e(L,H,K)e(L,H,K)^{g}=0,~{\rm~if~}g\not\in N_{G}(D).\end{equation} Let $g \not\in N_{G}(D)$. We have $$e(L,H,K)e(L,H,K)^{g}=(\sum_{x \in T}\varepsilon(H,K)^{x})(\sum_{y \in T}\varepsilon(H,K)^{y})^{g},$$ where $T$ is a transversal of $\operatorname{Cen}_{L}(\varepsilon(H,K))$ in $L.$ Thus eqn (\ref{e18}) follows if we show that for $$\varepsilon(H,K)^{x}\varepsilon(H,K)^{yg}=0, {\rm~~~~for~all~} x,~y \in T.$$ By Lemma \ref{l6}, $\varepsilon(H,K)\varepsilon(A,D)=\varepsilon(H,K)$. This gives $\varepsilon(H,K)^{x}\varepsilon(A,D)^{x}=\varepsilon(H,K)^{x}$ and $\varepsilon(A,D)^{yg}\varepsilon(H,K)^{yg}=\varepsilon(H,K)^{yg}.$  As $\psi_{A}$ is invariant in $L$, we have $D \unlhd L$ and hence $\varepsilon(A,D)^{x} = \varepsilon(A,D)$ and $\varepsilon(A,D)^{yg} = \varepsilon(A,D)^{g}$. Thus,  $\varepsilon(H,K)^{x}\varepsilon(H,K)^{yg}\linebreak = \varepsilon(H,K)^{x}\varepsilon(A,D)\varepsilon(A,D)^{g} \varepsilon(H,K)^{yg} = 0,$ as $g\not\in N_{G}(D)$. This proves eqn (\ref{e17}), which also yields \begin{equation}\label{e19}\operatorname{Cen}_{G}(e_{\mathbb{Q}}(\psi^{L}))\leqslant N_{G}(D).\end{equation} Now, let $g \in N_{G}(D)\setminus \operatorname{Cen}_{G}(e_{\mathbb{Q}}(\psi^{L})).$ Since $\psi^{L} \in \operatorname{Irr}L$ and $L\unlhd N_{G}(D)$, it follows that $e_{Q}(\psi^{L})$ and $e_{Q}(\psi^{L})^{g}$ are distinct primitive central idempotents of $\mathbb{Q}L$ and therefore $e_{\mathbb{Q}}(\psi^{L})e_{\mathbb{Q}}(\psi^{L})^{g}= 0$. This proves (i).\vspace{.2cm}\\(ii) Let $T$, $T_{1}$, $T_{2}$, $T_{3}$ respectively be a right transversal of $\operatorname{Cen}_{L}(\varepsilon(H,K))$ in $G$, $\operatorname{Cen}_{L}(\varepsilon(H,K))$ in $L$, $L$ in $\operatorname{Cen}_{G}(e_{\mathbb{Q}}(\psi^{L}))$, $\operatorname{Cen}_{G}(e_{\mathbb{Q}}(\psi^{L}))$ in $G$. We have, $$\begin{array}{lll}\displaystyle\sum_{g \in T}\varepsilon(H,K)^{g}&=&\displaystyle\sum_{z \in T_{3}}\sum_{y \in T_{2}}(\sum_{x \in T_{1}}\varepsilon(H,K)^{x})^{yz}\\&=&\displaystyle\sum_{z \in T_{3}}\sum_{y \in T_{2}}e(L,H,K)^{yz}\\&=&\displaystyle\frac{[\operatorname{Cen}_{G}(e_{\mathbb{Q}}(\psi^{L})):L]}{\alpha_{(L,H,K)}}\sum_{z \in T_{3}}e_{\mathbb{Q}}(\psi^{L})^{z}.\end{array} $$ Also, it is easy to see that $$ \begin{array}{lll} \displaystyle\sum_{g \in T}\varepsilon(H,K)^{g}&=&[\operatorname{Cen}_{G}(\varepsilon(H,K)):\operatorname{Cen}_{L}(\varepsilon(H,K))]e(G,H,K)\\ &=& \frac{[\operatorname{Cen}_{G}(\varepsilon(H,K)):\operatorname{Cen}_{L}(\varepsilon(H,K))]}{\alpha_{(G,H,K)}}e_{\mathbb{Q}}(\psi^{G}).\end{array}$$  We know from part(i) that $\displaystyle\sum_{z \in T_{3}}e_{\mathbb{Q}}(\psi^{L})^{z}$ is an idempotent. Also $e_{\mathbb{Q}}(\psi^{G})$ being an idempotent, it follows immediately that  $\frac{[\operatorname{Cen}_{G}(e_{\mathbb{Q}}(\psi^{L})):L]}{\alpha_{(L,H,K)}}= \frac{[\operatorname{Cen}_{G}(\varepsilon(H,K)):\operatorname{Cen}_{L}(\varepsilon(H,K))]}{\alpha_{(G,H,K)}}$, which gives the desired result.
 $~\Box$ \vspace{.3cm} \par The above proposition gives the following generalization of (\cite{OdRS04}, Corollary 3.6): \begin{cor}\label{c1} Let $(H,K)$ be a pair of subgroups of a finite group $G$ and $A$ be a normal subgroup of $G$ contained in $H$ satisfying the following conditions: \begin{description}\item[(i)] $K\unlhd H\unlhd N_{G}(D)$, where $D=K\cap A$;\item[(ii)] $H/K$ is cyclic and a maximal abelian subgroup of $N_{G}(K)/K$.\end{description} Then $(H,K)$ is a strong Shoda pair of $G$ and $e(G,H,K)$ is a primitive central idempotent of $\mathbb{Q}G.$\end{cor}\noindent{\bf Proof.} In view of (i), $H=L$ satisfies the hypothesis of the above proposition. Therefore, the above proposition yields that the distinct $G$-conjugates of $\varepsilon(H,K)$ are mutually orthogonal and hence $(H,K)$ is a strong Shoda pair of $G$. $~\Box$ \vspace{.3cm} \par  We also have the following: \begin{cor}\label{c2} Let $G \in \mathcal{C}$, $N$ a normal subgroup of $G$. Let $(H,H,\vartheta) \in \mathcal{L}_{N}$. If  $H\unlhd N_{G}(\operatorname{ker}\vartheta_{2})$,  then $(H,\operatorname{ker}\vartheta)$ is a strong Shoda pair of $G$, where $(H_{2}, A_{2}, \vartheta_{2})$ is the second node of $(H,H,\vartheta)$.  \end{cor}\noindent{\bf Proof.} Let $K=\operatorname{ker}\vartheta.$ From Theorem \ref{t2}, $(H,K)$ is a Shoda pair of $G$. By considering $A=A_{2}$, $L=H$ and $\psi=\vartheta$ in the above proposition, it follows that the distinct $G$-conjugates of $\varepsilon(H,K)$ are mutually orthogonal. Also $\vartheta_{A_{2}}$=$\vartheta_{2}$ implies that $K\cap A_{2}=\operatorname{ker}\vartheta_{2}$ and hence $N_{G}(K)\leqslant N_{G}(\operatorname{ker}\vartheta_{2})$. Now $H$ being normal in $N_{G}(\operatorname{ker}\vartheta_{2})$, it follows that $H\unlhd N_{G}(K)$. Consequently, by (\cite{OdRS04}, Proposition 3.3), it follows that $(H,K)$ is a strong Shoda pair of $G$. $~\Box$ \begin{remark}\label{r2} From the above corollary, it follows immediately that if $G \in \mathcal{C}$ and $(H,H,\vartheta) \in \bigcup_{N \in \mathcal{N}}\mathcal{L}_{N}$ is of height $2$, then $(H, \operatorname{ker}\vartheta)$ is a strong Shoda pair of $G$.\end{remark}
\noindent{\bf Proof of Theorem \ref{t4}} (i) By taking $L=A=H_{n}$ and $G=H_{n-1}$ in Proposition \ref{p1}, it follows that $$\alpha_{(H_{n-1},H,K)}=1.$$ Also for $2\leq i\leq n-1$, by taking $L=H_{i}$, $G=H_{i-1}$, $A=A_{i}$, we have  $$\alpha_{(H_{i-1},H,K)}=\alpha_{(H_{i},H,K)}\frac{[\operatorname{Cen}_{H_{i-1}}(\varepsilon(H,K)):\operatorname{Cen}_{H_{i}}(\varepsilon(H,K))]}{
 [\operatorname{Cen}_{H_{i-1}}(e(H_{i},H,K)):H_{i}]}.$$ Consequently
 $$ \begin{array}{lll} \alpha_{(G,H,K)} &= & \prod_{2\leq i\leq n-1}\frac{[\operatorname{Cen}_{H_{i-1}}(\varepsilon(H,K)):\operatorname{Cen}_{H_{i}}(\varepsilon(H,K))]}
 {[\operatorname{Cen}_{H_{i-1}}(e(H_{i},H,K)):H_{i}]}\\ & = & \frac{[\operatorname{Cen}_{G}(\varepsilon(H,K)):\operatorname{Cen}_{H_{n-1}}(\varepsilon(H,K))]}{\prod_{2\leq i\leq n-1}[\operatorname{Cen}_{H_{i-1}}(e(H_{i},H,K)):H_{i}]}, \end{array}  $$ as desired.
 \\ (ii) In view of (i), $\alpha_{(G,H,K)}=1$ if, and only if, $$\displaystyle\prod_{2\leq i\leq n-1}[\operatorname{Cen}_{H_{i-1}}(\varepsilon(H,K)):\operatorname{Cen}_{H_{i}}(\varepsilon(H,K))]=\displaystyle\prod_{2\leq i\leq n-1}[\operatorname{Cen}_{H_{i-1}}(e(H_{i},H,K)):H_{i}].$$ But $\operatorname{Cen}_{H_{i-1}}(\varepsilon(H,K))/\operatorname{Cen}_{H_{i}}(\varepsilon(H,K))$ being isomorphic to a subgroup of \linebreak $\operatorname{Cen}_{H_{i-1}}(e(H_{i},H,K))/H_{i},$ the above equation holds if, and only if, $$[\operatorname{Cen}_{H_{i-1}}(\varepsilon(H,K)):\operatorname{Cen}_{H_{i}}(\varepsilon(H,K))]=[\operatorname{Cen}_{H_{i-1}}(e(H_{i},H,K)):H_{i}],$$ for all i, $2\leq i\leq n-1,$ which yields the required result. \\(iii) Let $2\leq i \leq n-1$. Clearly $x  \in \operatorname{Cen}_{H_{i-1}}(e(H_{i},H,K))$ if, and only if, $ x \in \operatorname{Cen}_{H_{i-1}}(e_{\mathbb{Q}}(\vartheta^{H_{i}})) $ if, and only if,  $e_{\mathbb{Q}}((\vartheta^{H_{i}})^{x})=e_{\mathbb{Q}}(\vartheta^{H_{i}})$ if, and only if, $(\vartheta^{H_{i}})^{x}= \sigma\circ \vartheta^{H_{i}}$ for some $\sigma \in \operatorname{\mathcal{A}ut}(\mathbb{C})$. However, the later holds if, and only if, $x \in N_{H_{i-1}}(\operatorname{ker}\vartheta_{i})$, provided $(H, H, \vartheta)$ is good. This proves (iii).
 $~\Box$ \section{Simple components} Given $G \in \mathcal{C}$, it follows from Theorem \ref{t2} that any simple component of  $\mathbb{Q}G$ is given by $ \mathbb{Q}G e_{\mathbb{Q}}(\vartheta^{G})$, where $(H,H,\vartheta) \in \cup_{N \in \mathcal{N}}\mathcal{L}_{N}$. Let's now compute the structure of $\mathbb{Q}G e_{\mathbb{Q}}(\vartheta^{G})$.\par For a ring $R$, let $\mathcal{U}(R)$ be the unit group of $R$ and $M_{n}(R)$ the ring of $n \times n$ matrices over $R$. Denote by $R*_{\tau}^{\sigma}G$, the crossed product of the group $G$ over the ring $R$ with action $\sigma$ and twisting $\tau$.\vspace{.2cm}\\
   We begin with the following generalization of Proposition 3.4 of \cite{OdRS04}: \begin{prop}\label{p2} Let $G$ be a finite group and $(H,K)$  a Shoda pair of $G$. Let $\psi,~A,~L$ be as in Proposition \ref{p1}. Then  $$\mathbb{Q}Ge_{\mathbb{Q}}(\psi^{G})\cong M_{n}(\mathbb{Q}Le_{\mathbb{Q}}(\psi^{L})\ast^{\sigma}_{\tau}\operatorname{Cen}_{G}(e_{\mathbb{Q}}(\psi^{L}))/L),$$ where $n=[G:\operatorname{Cen}_{G}(e_{\mathbb{Q}}(\psi^{L}))]$, the action $\sigma$ {\rm:}$\operatorname{Cen}_{G}(e_{\mathbb{Q}}(\psi^{L}))/L\rightarrow\operatorname{\mathcal{A}ut}(\mathbb{Q}Le_{\mathbb{Q}}(\psi^{L}))$ maps $\overline{x}$ to the conjugation automorphism induced by $x$ and the twisting \linebreak $\tau${\rm:}$\operatorname{Cen}_{G}(e_{\mathbb{Q}}(\psi^{L}))/L\times \operatorname{Cen}_{G}(e_{\mathbb{Q}}(\psi^{L}))/L\rightarrow
\mathcal{U}(\mathbb{Q}Le_{\mathbb{Q}}(\psi^{L}))$ is given by $(\overline{g_{1}},\overline{g_{2}})\mapsto g$, where $g \in L$ is such that $\overline{g_{1}} \cdot\overline{g_{2}}=g\cdot\overline{g_{1}}\overline{g_{2}}$ for $\overline{g_{1}},\overline{g_{2}} \in \operatorname{Cen}_{G}(e_{\mathbb{Q}}(\psi^{L}))/L.$ \end{prop}\noindent{\bf Proof.} Clearly $\mathbb{Q}Ge_{\mathbb{Q}}(\psi^{G})$ is isomorphic to the ring  $\operatorname{End}_{\mathbb{Q}G}(\mathbb{Q}Ge_{\mathbb{Q}}(\psi^{G}))$ of $\mathbb{Q}G$ endomorphisms of $\mathbb{Q}Ge_{\mathbb{Q}}(\psi^{G})$. As $e_{\mathbb{Q}}(\psi^{G})$ is the sum of distinct $G$-conjugates of $e_{\mathbb{Q}}(\psi^{L})$, we have $$\operatorname{End}_{\mathbb{Q}G}(\mathbb{Q}Ge_{\mathbb{Q}}(\psi^{G}))\cong M_{n}(\operatorname{End}_{\mathbb{Q}G}(\mathbb{Q}Ge_{\mathbb{Q}}(\psi^{L}))),$$ where $n=[G:\operatorname{Cen}_{G}(e_{\mathbb{Q}}(\psi^{L}))]$. Also the map $f \mapsto f(e_{\mathbb{Q}}(\psi^{L}))$ yields isomorphism $$\operatorname{End}_{\mathbb{Q}G}(\mathbb{Q}Ge_{\mathbb{Q}}(\psi^{L}))\cong e_{\mathbb{Q}}(\psi^{L})\mathbb{Q}G e_{\mathbb{Q}}(\psi^{L}).$$ But distinct $G$-conjugates of $e_{\mathbb{Q}}(\psi^{L})$ being orthogonal, we have $e_{\mathbb{Q}}(\psi^{L})\mathbb{Q}G e_{\mathbb{Q}}(\psi^{L})\cong e_{\mathbb{Q}}(\psi^{L})\mathbb{Q}C e_{\mathbb{Q}}(\psi^{L})=\mathbb{Q}C e_{\mathbb{Q}}(\psi^{L}),$ where $C=\operatorname{Cen}_{G}(e_{\mathbb{Q}}(\psi^{L})).$ Consequently, $\mathbb{Q}Ge_{\mathbb{Q}}(\psi^{G})$ is isomorphic to $ M_{n}(\mathbb{Q}Ce_{\mathbb{Q}}(\psi^{L})).$ Since $L\unlhd N_{G}(\operatorname{ker}\psi_{A})$, eqn (\ref{e19}) gives $L\unlhd C$ and hence $$\mathbb{Q}Ce_{\mathbb{Q}}(\psi^{L})\cong \mathbb{Q}Le_{\mathbb{Q}}(\psi^{L})\ast^{\sigma}_{\tau}C/L,$$ where $\sigma$ and $\tau$ are as in the statement. This completes the proof.
 $~\Box$ \vspace{.3cm} \par We now proceed to compute the structure of the simple component $\mathbb{Q}G e_{\mathbb{Q}}(\vartheta^{G})$ of $\mathbb{Q}G$, where $G \in \mathcal{C}$ and $(H,H,\vartheta) \in \cup_{N \in \mathcal{N}}\mathcal{L}_{N}$. \par Suppose that $(H,H,\vartheta) \in \cup_{N \in \mathcal{N}}\mathcal{L}_{N}$ is of height $n$ and  $(H_{i},A_{i},\vartheta_{i})$ is its $i^{th}$ node, $1\leq i\leq n$. Let  $K = \operatorname{ker}\vartheta$. Let $1\leq i\leq n-1$. We must notice that $G=H_{i}$, $L=H_{i+1}$, $A=A_{i+1}$ and $\psi=\vartheta$ satisfy the hypothesis of Proposition \ref{p2}. Let $k_{i}=[H_{i}:\operatorname{Cen}_{H_{i}}(e_{\mathbb{Q}}(\vartheta^{H_{i+1}}))],$ $\sigma_{i}$ the action of $\operatorname{Cen}_{H_{i}}(e_{\mathbb{Q}}(\vartheta^{H_{i+1}}))/H_{i+1}$ on $\mathbb{Q}H_{i+1}e_{\mathbb{Q}}(\vartheta^{H_{i+1}})$ by conjugation, $\tau_{i}$: $\operatorname{Cen}_{H_{i}}(e_{\mathbb{Q}}(\vartheta^{H_{i+1}}))/H_{i+1}\times \operatorname{Cen}_{H_{i}}(e_{\mathbb{Q}}(\vartheta^{H_{i+1}}))/H_{i+1}\rightarrow \mathcal{U}(\mathbb{Q}H_{i+1}e_{\mathbb{Q}}(\vartheta^{H_{i+1}}))$  the twisting given by $(\overline{g_{1}},\overline{g_{2}})\mapsto g$, where $g \in H_{i+1}$ is such that $\overline{g_{1}} \cdot\overline{g_{2}}=g\cdot\overline{g_{1}}\overline{g_{2}}$, for $\overline{g_{1}},\overline{g_{2}} \in \operatorname{Cen}_{H_{i}}(e_{\mathbb{Q}}(\vartheta^{H_{i+1}}))/H_{i+1}.$ Observe that $\operatorname{Cen}_{H_{i}}(e_{\mathbb{Q}}(\vartheta^{H_{i+1}}))$=$\operatorname{Cen}_{H_{i}}(e(H_{i+1},H,K))$ for all $i$. \par Apply Proposition \ref{p2} with $G=H_{n-1}$ and $A=L=H_{n}$. It follows that  \begin{equation} \label{e20}  \mathbb{Q}H_{n-1}e_{\mathbb{Q}}(\vartheta^{H_{n-1}}) \cong M_{k_{n-1}}(\mathbb{Q}H\varepsilon(H,K)\ast^{\sigma_{n-1}}_{\tau_{n-1}} N_{H_{n-1}}(K)/H),\end{equation} as $ e_{\mathbb{Q}}(\vartheta^{H_{n}}) = \varepsilon(H,K)$ and $\operatorname{Cen}_{H_{n-1}}(\varepsilon(H,K))= N_{H_{n-1}}(K)$. Denote by $R$ the matrix ring on right hand side of the above isomorphism. \par The action  of $\operatorname{Cen}_{H_{n-2}}(e(H_{n-1},H,K))/H_{n-1}$ on $\mathbb{Q}H_{n-1}e_{\mathbb{Q}}(\vartheta^{H_{n-1}})$  given by $\sigma_{n-2}$ induces its natural action  on $R$ given by $\overline{x}\mapsto \eta\circ\sigma_{n-2}(\overline{x})\circ\eta^{-1}$, where $\overline{x}\in \operatorname{Cen}_{H_{n-2}}(e(H_{n-1},H,K))/H_{n-1}$ and $\eta$ is the isomorphism in eqn (\ref{e20}). For notational convenience, we denote this action on $R$ again by $\sigma_{n-2}$. Similarly, the twisting $\eta\circ\tau_{n-2}$ from  $\operatorname{Cen}_{H_{n-2}}(e(H_{n-1},H,K))/H_{n-1}\times \operatorname{Cen}_{H_{n-2}}(e(H_{n-1},H,K))/H_{n-1}$ to $\mathcal{U}(R)$ will again be denoted by $\tau_{n-2}$ for convenience. Applying Proposition \ref{p2} with $G=H_{n-2}$, $L=H_{n-1}$, $A = A_{n-1}$, it follows that  $\mathbb{Q}H_{n-2}e_{\mathbb{Q}}(\vartheta^{H_{n-2}})$ is isomorphic to $$M_{k_{n-2}}((M_{k_{n-1}}(\mathbb{Q}H\varepsilon(H,K)\ast^{\sigma_{n-1}}_{\tau_{n-1}} N_{H_{n-1}}(K)/H)\ast^{\sigma_{n-2}}_{\tau_{n-2}}\operatorname{Cen}_{H_{n-2}}(e(H_{n-1},H,K)/H_{n-1}).$$ \par Continue this process, after $n-1$ steps, we obtain the following:  \begin{theorem}\label{t5} Let $G \in \mathcal{C}$ and $(H,H,\vartheta) \in \cup_{N \in \mathcal{N}}\mathcal{L}_{N}$. If   $(H,H,\vartheta)$ is of height  $n$ and $(H_{i},A_{i},\vartheta_{i})$ is  its $i^{th}$ node, $1\leq i\leq n$, then  $\mathbb{Q}Ge_{\mathbb{Q}}(\psi^{G})$ is isomorphic to  $$M_{k_{1}}(M_{k_{2}}\cdots(M_{k_{n-1}}(\mathbb{Q}(\xi_{k})\ast^{\sigma_{n-1}}_{\tau_{n-1}}N_{H_{n-1}}(K)/H)\ast^{\sigma_{n-2}}_{\tau_{n-2}}
\cdots\ast^{\sigma_{1}}_{\tau_{1}}\operatorname{Cen}_{H_{1}}(e(H_{2},H,K))/H_{2}),$$ where $k$=$[H:K]$, $\xi_{k}$ is a primitive $k^{th}$ root of unity, and $\sigma_{i},~\tau_{i},~k_{i}$ are as defined above.\end{theorem}
\section{Examples} In this section, we illustrate the construction of Shoda pairs. We begin by observing the following facts for an arbitrary group $G$ in $\mathcal{C}$: \begin{enumerate} \item The directed tree $\mathcal{G}_{N}$, when $N=G$, is just the vertex $(G,G,1_{G})$, as $Cl(G,G,1_{G})$ is an empty set. In this case,   $\mathcal{G}_{N}$  corresponds to the Shoda pair $(G,G)$. \item If $N$ is a normal subgroup of $G$ such that $\mathcal{A}_{(G, N,1_{N})}/N$ is cyclic,  then \begin{quote} $ (G, N ,1_{N})$ has only one d.c.c., namely, $(I_{G}(\varphi),\mathcal{A}_{(G,N,1_{N})}, \varphi)$, \end{quote} where $\varphi$ can be taken to be any linear character on $\mathcal{A}_{(G,N,1_{N})}$ with kernel $N$. To show this, consider $\varphi \in \operatorname{\widetilde{Lin}}(\mathcal{A}_{(G,N,1_{N})}|1_{N})$. We have   $\varphi_{N} = 1_{N}$  and $\operatorname{ker}\varphi^{G} = N.$  The condition $\varphi_{N} = 1_{N}$ gives $N \leq \operatorname{ker} \varphi \leq  \mathcal{A}_{(G,N,1_{N})}$. This yields that $\operatorname{ker}\varphi$ is a normal subgroup of $G$, as $\mathcal{A}_{(G, N, 1_{N})}/N$  is a cyclic normal subgroup of $G/N$. Consequently, we have  $\operatorname{ker} \varphi = \operatorname{ker}\varphi^{G} = N$, and thus  $\varphi$ is a linear character on $\mathcal{A}_{(G,N,1_{N})}$ with kernel $N$. Conversely, it is clear that any linear character on $\mathcal{A}_{(G,N,1_{N})}$ with kernel $N$ belongs to $\operatorname{\widetilde{Lin}}(\mathcal{A}_{(G,N,1_{N})}|1_{N})$. Hence $\operatorname{\widetilde{Lin}}(\mathcal{A}_{(G, N, 1_{N})}|1_{N})$ consists precisely  of all  linear characters on $\mathcal{A}_{(G,N,1_{N})}$ with kernel $N$.  It is easy to see that all these characters lie in the same orbit under the double action. Consequently, there is only one d.c.c. of $(G,N,1_{N})$, as desired. \par  \indent Further, if $\mathcal{A}_{(G,N,1_{N})}$ is such that $\mathcal{A}_{(G,N,1_{N})}/N$ is maximal among all the abelian subgroups of $G/N$, then  we  have $I_{G}(\varphi)= \mathcal{A}_{(G,N,1_{N})}$ and hence there is no further d.c.c. of $(I_{G}(\varphi),\mathcal{A}_{(G,N,1_{N})}, \varphi)$ and the process stops here and the corresponding directed tree is as follows: \begin{figure}[H] \centering \begin{tikzpicture}[font=\footnotesize, edge from parent/.style={draw,thick}]
  \node{$(G,N,1_{N})$}[grow'=up]
     child {node {$(\mathcal{A}_{(G,N,1_{N})}, \mathcal{A}_{(G,N,1_{N})},\varphi)$~~~~}edge from parent[-stealth]}
     ;
\end{tikzpicture} \caption{$\mathcal{G}_{N}$} \end{figure} If $\mathcal{A}_{(G,N,1_{N})}/N$ is not maximal among all the abelian subgroups of $G/N$, then $I_{G}(\varphi) \neq \mathcal{A}_{(G,N,1_{N})}$. In this case, we further need to compute the d.c.c. of $(I_{G}(\varphi),\mathcal{A}_{(G,N,1_{N})}, \varphi)$ in order to determine $\mathcal{G}_{N}$.\item From the above fact, it follows that if $N$ is a normal subgroup of $G$ with $G/N$ cyclic, then its corresponding directed tree is as follows:
\begin{figure}[H] \centering \begin{tikzpicture}[font=\footnotesize, edge from parent/.style={draw,thick}]
  \node{$(G,N,1_{N})$}[grow'=up]
     child {node {$(G,G,\varphi)$~~~~}edge from parent[-stealth]}
     ;
\end{tikzpicture} \caption{$\mathcal{G}_{N}$} \end{figure}
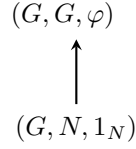
where $\varphi$ is any linear character of $G$  with $\operatorname{ker}\varphi$ = $N$ and it yields the Shoda pair $(G,N)$ of $G$.\item  If $N$ is a normal subgroup of $G$  with $G/N$ abelian but not cyclic then there is no d.c.c. of $(G,N,1_{N})$. As, in this case $\mathcal{A}_{(G,N,1_{N})}=G$ and $\varphi \in \operatorname{\widetilde{Lin}}(G|1_{N})$ implies that $\varphi_{N} = 1_{N}$  and $\operatorname{ker}\varphi^{G} = N$ gives $\operatorname{ker}\varphi = N$. Consequently, $G/N$ is cyclic, which is not the case. Hence the directed tree $\mathcal{G}_{N}$ is just the vertex $(G,N,1_{N})$, which does not yield any Shoda pair as there is no leaf of the required type.

     \item Given $N$-linear character triple $(H, A, \vartheta)$ of $G$, consider the set $\mathcal{K}$ of all normal subgroups $K$ of $\mathcal{A}_{(H, A, \vartheta)}$ satisfying the following: \begin{description} \item(i) $\mathcal{A}_{(H, A, \vartheta)}/K$ is cyclic; \item(ii) $K\cap A= \operatorname{ker}\vartheta$; \item(iii) $ \operatorname{core}_{G}(K) = N.$  \end{description} Let $K_{1}$, $K_{2}$, $\cdots$, $K_{r}$ be a set of representatives of $\mathcal{K}$ under the equivalence relation defined by conjugacy  of subgroups in $H$. If $\varphi_{i}$ is a linear character on $\mathcal{A}_{(H, A, \vartheta)}$ with kernel $K_{i}$, $1 \leq i \leq r$, then $$Cl{(H,A,\vartheta)} = \{ (I_{H}(\varphi_{i}), \mathcal{A}_{(H, A, \vartheta)}, \varphi_{i}) ~|~ 1 \leq i \leq r\}.$$ \noindent To show this, consider $\varphi \in  \operatorname{\widetilde{Lin}}(\mathcal{A}_{(H, A,  \vartheta )}| \vartheta)$. Let $K = \operatorname{ker}  \varphi$. Clearly $\mathcal{A}_{(H, A,  \vartheta )}/K$ is cyclic. Also $ \varphi_{A}= \vartheta$ and $ \operatorname{ker}  \varphi^{G}=N$ implies that $K\cap A=  \operatorname{ker}  \vartheta$ and $\operatorname{core}_{G}(K)=N$. Hence $K \in  \mathcal{K}$. Thus we have shown that if  $\varphi \in  \operatorname{\widetilde{Lin}}(\mathcal{A}_{(H, A,  \vartheta )}| \vartheta)$, then $\operatorname{ker} \varphi \in \mathcal{K}$. Conversely, it is clear that if $K \in \mathcal{K}$, then any linear character on $\mathcal{A}_{(H, A,  \vartheta )}$ with kernel $K$ lie in $\operatorname{\widetilde{Lin}}(\mathcal{A}_{(H, A,  \vartheta )}| \vartheta)$. Furthermore, note that  $\varphi_{1}$ and $\varphi_{2} \in \operatorname{\widetilde{Lin}}(\mathcal{A}_{(H, A,  \vartheta )}| \vartheta)$ lie in the same orbit under the double action if, and only if,  $\operatorname{ker}  \varphi_{1}$ and  $\operatorname{ker}  \varphi_{2}$  are conjugate in $H$. This  yields the desired result. \item In the above fact, if $(H, A, \vartheta)$ is such that $H/\operatorname{ker} \vartheta $  is abelian, then it may be noted that $K$, $K' \in \mathcal{K}$ are conjugate in $H$ if, and only if, $K=K'$. \end{enumerate}

    \subsection{Example $1$} Let $G$  be the group generated by $x_{i},~1\leq i\leq 6$, with the following defining relations:
\begin{quote}  $x_{1}^{2}x_{2}^{-1}$=$x_{2}^{2}x_{3}^{-1}$=$x_{4}^{5}$=$x_{3}^{2}$=$x_{5}^{5}$=$x_{6}^{5}$=$1,$\\
               $[x_{2},x_{1}]$=$[x_{3},x_{1}]$=$[x_{3},x_{2}]$=$[x_{6},x_{3}]$=$[x_{6},x_{4}]$=$[x_{6},x_{5}]$=$1,$\\
                $[x_{5},x_{4}]$=$x_{6}$, $[x_{5},x_{1}]$=$x_{4}x_{5},$\\
                $[x_{6},x_{1}]$=$x_{6}^{2}$, $[x_{4},x_{2}]$=$x_{4}x_{6}^{2}$, $[x_{6},x_{2}]$=$x_{6}^{3},$\\
                 $[x_{5},x_{2}]$=$x_{5}x_{6}^{2}$, $[x_{5},x_{3}]$=$x_{5}^{3}x_{6}^{2}$, $[x_{4},x_{3}]$=$x_{4}^{3}x_{6}^{2}$, $[x_{4},x_{1}]$=$x_{4}^{2}x_{5}^{3}x_{6}^{4}.$\end{quote}
                 This group is SmallGroup(1000,86) in GAP library. We have already mentioned in section $2$ that it belongs to $\mathcal{C}$ but it is not strongly monomial. There are $6$ normal subgroups of $G$ given by $N_{1}$=$G$, $N_{2}$=$\langle x_{2},x_{3},x_{4},x_{5},x_{6}\rangle$, $N_{3}$=$\langle x_{3},x_{4},x_{5},x_{6}\rangle$, $N_{4}$=$\langle x_{4},x_{5},x_{6}\rangle$, $N_{5}$=$\langle x_{6}\rangle$, $N_{6}$=$\langle1\rangle.$ We will compute the directed tree $\mathcal{G}_{N_{i}}$, for all $i$. From fact $1$, $\mathcal{G}_{N_{1}}$ is just the vertex $(G,G,1_{G})$ and it corresponds to the Shoda pair $(G,G)$. For $2\leq i\leq 4$, $G/N_{i}$ is cyclic. Therefore if $\varphi_{1},~\varphi_{2},~\varphi_{3}$ are linear characters on $G$ with kernel $N_{2},~N_{3},~N_{4}$ respectively, then, by fact $3$ above, we have the following:
               \begin{figure}[H]
\centering
\begin{minipage}{.3\textwidth}\centering
 \begin{tikzpicture}[font=\footnotesize, edge from parent/.style={draw,thick}]
\node{$(G,N_{2},1_{N_{2}})$}[grow'=up]
  child {node{$(G,G,\varphi_{1})$} edge from parent[-stealth]}
  ;
\end{tikzpicture}
 \captionof{figure}{$\mathcal{G}_{N_{2}}$}
\end{minipage}%
\begin{minipage}{.3\textwidth}
  \centering
\begin{tikzpicture}[font=\footnotesize, edge from parent/.style={draw,thick}]
\node{$(G,N_{3},1_{N_{3}})$}[grow'=up]
  child {node{$(G,G,\varphi_{2})$} edge from parent[-stealth]}
  ;
\end{tikzpicture}
\captionof{figure}{$\mathcal{G}_{N_{3}}$}
\end{minipage}%
\begin{minipage}{.3\textwidth}
  \centering
\begin{tikzpicture}[font=\footnotesize, edge from parent/.style={draw,thick}]
\node{$(G,N_{4},1_{N_{4}})$}[grow'=up]
  child {node{$(G,G,\varphi_{3})$} edge from parent[-stealth]}
  ;
\end{tikzpicture}
\captionof{figure}{$\mathcal{G}_{N_{4}}$}
\end{minipage}
\end{figure}
       \noindent and the leaves $(G,G,\varphi_{1})$, $(G,G,\varphi_{2})$ and $(G,G,\varphi_{3})$   yield the Shoda pairs $(G, N_{2})$, $(G, N_{3})$, $(G, N_{4})$ of $G$ respectively.  \vspace{.2cm}\\
 \noindent $\underline{{\rm Construction~of~}\mathcal{G}_{N_{5}}:}$ We first need to compute d.c.c. of $(G,N_{5},1_{N_{5}})$. Observe that $N_{4}/N_{5}$ is an abelian normal subgroup of maximal order in $G/N_{5}$. Therefore, we set $\mathcal{A}_{(G,N_{5},1_{N_{5}})}= N_{4}$. We now use fact 5 to compute $Cl(G, N_{5}, 1_{N_{5}})$.  It turns out that $\mathcal{K} = \{ \langle x_{4}^{i}x_{5}, x_{6}\rangle~|~ 0\leq i \leq 4 \}\cup\{\langle x_{4},x_{6}\rangle\}$. Moreover,  $\langle x_{5},x_{6}\rangle$, $\langle x_{4},x_{6}\rangle$, $\langle x_{4}^{-1}x_{5},x_{6}\rangle$ are the only  subgroups in $\mathcal{K}$ which are distinct up to conjugacy in $G$. Consider the linear characters $\varphi_{1}$, $\varphi_{2}$, $\varphi_{3}$ on $N_{4}$ given as follows: \begin{center} $\varphi_{1} : x_{4}\mapsto \xi_{5},~x_{5} \mapsto 1,~x_{6} \mapsto 1$\\
                                                         $\varphi_{2} : x_{4} \mapsto 1,~x_{5} \mapsto \xi_{5},~x_{6} \mapsto 1$\\
                                                         $\varphi_{3} : x_{4} \mapsto \xi_{5},~x_{5} \mapsto \xi_{5},~x_{6} \mapsto 1$,\end{center} where $\xi_{5}$ is a primitive $5^{th}$ root of unity. Clearly the kernels of $\varphi_{1}$, $\varphi_{2}$ and $\varphi_{3}$ are
 $\langle x_{5},x_{6}\rangle$, $\langle x_{4},x_{6}\rangle$ and $\langle x_{4}^{-1}x_{5},x_{6}\rangle$ respectively. Also we have $
 I_{G}(\varphi_{1})=I_{G}(\varphi_{2})=I_{G}(\varphi_{3})=N_{4}.$  Hence, the directed tree is as follows:
 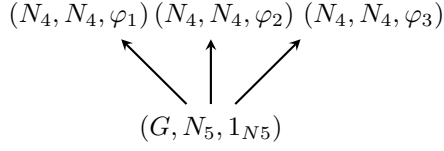
\begin{figure}[H] \centering \begin{tikzpicture}[font=\footnotesize, edge from parent/.style={draw,thick}]
  \node{$(G,N_{5},1_{N5})$}[grow'=up]
     child {node {$(N_{4},N_{4},\varphi_{1})$~~~~~}edge from parent[-stealth]}
     child {node {~~~~~~~$(N_{4},N_{4},\varphi_{2})$~~~~~}edge from parent[-stealth]}
     child {node {~~~~~~~~~~~$(N_{4},N_{4},\varphi_{3})$}edge from parent[-stealth]}
     ;
\end{tikzpicture} \caption{$\mathcal{G}_{N_{5}}$} \end{figure}
  \noindent The three leaves $(N_{4},N_{4},\varphi_{1})$, $(N_{4},N_{4},\varphi_{2})$, $(N_{4},N_{4},\varphi_{3})$, respectively, yield  Shoda pairs $(\langle x_{4},x_{5},x_{6}\rangle,\langle x_{5},x_{6}\rangle) $, $(\langle x_{4},x_{5},x_{6}\rangle,\langle x_{4},x_{6}\rangle)$ and $(\langle x_{4},x_{5},x_{6}\rangle,\langle x_{4}^{-1}x_{5},x_{6}\rangle)$ of $G$. \vspace{.2cm}\\ $\underline{{\rm Construction~of~}\mathcal{G}_{N_{6}}:}$ Recall that $N_{6} =\langle 1 \rangle$. Note that $N_{5}$ is an abelian normal subgroup of maximal order in $G$. Therefore, we  set  $\mathcal{A}_{(G,N_{6},1_{N_{6}})} = N_{5}$. As  $N_{5}/N_{6}$ is cyclic, by fact $2$, $$Cl(G, N_{6}, 1_{N_{6}}) = \{ (G_{1}, N_{5}, \varphi_{1})\},$$ where $\varphi_{1}:N_{5} \rightarrow \mathbb{C} $ maps $x_{6}$  to $\xi_{5}$, and $G_{1}= I_{G}(\varphi_{1})= \langle x_{3},x_{4},x_{5},x_{6}\rangle. $
As $G_{1} \neq N_{5}$, we further need to compute d.c.c. of $(G_{1}, N_{5},\varphi_{1})$.  It is observed that $\langle x_{5},x_{6}\rangle/\operatorname{ker}\varphi_{1}$ is an abelian normal subgroup of maximal order in $G_{1}/\operatorname{ker}\varphi_{1}$. Set $\mathcal{A}_{(G_{1}, N_{5},\varphi_{1})}= \langle x_{5},x_{6}\rangle$. We now use fact $5$ to compute $Cl(G_{1}, N_{5},\varphi_{1})$. It turns out that there are $5$ subgroups in $\mathcal{K}$, namely $\langle x_{5}x_{6}^{i}\rangle$, $ 0\leq i \leq 4$, and all of them are conjugate in $G_{1}.$  Consider the linear character $\varphi_{2}$ on $\langle x_{5},x_{6}\rangle$ which maps $x_{5}$ to $1$ and $x_{6}$ to $\xi_{5}$. Hence $$Cl(G_{1}, N_{5},\varphi_{1}) = \{ (G_{2}, \langle x_{5},x_{6}\rangle , \varphi_{2})\},$$ where  $G_{2}= I_{G_{1}}(\varphi_{2})= \langle x_{5},x_{6},x_{3}x_{4}^{2}\rangle. $ Again $G_{2} \neq \langle x_{5},x_{6}\rangle $, and we compute $Cl(G_{2}, \langle x_{5},x_{6}\rangle , \varphi_{2})$. Now  $G_{2}/\operatorname{ker} \varphi_{2}$  is abelian. Therefore, using facts $5$ and $6$, we obtain that  $$Cl(G_{2}, \langle x_{5},x_{6}\rangle,\varphi_{2}) = \{ (G_{2}, G_{2}, \varphi_{3}), (G_{2}, G_{2}, \varphi_{4})\},$$  where $\varphi_{3}$ and $\varphi_{4}$ are linear characters on $G_{2}$ with kernel $\langle x_{3}x_{4}^{2}x_{6}^{3},x_{5}^{3}\rangle$ and $\langle x_{5}^{3}\rangle$ respectively. The process stops here and the corresponding tree is as follows: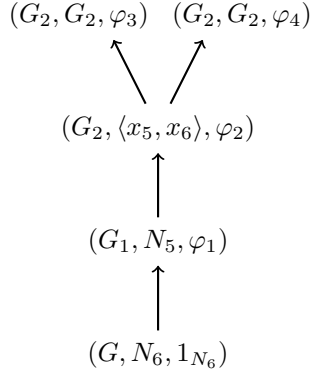
\begin{figure}[H] \centering \begin{tikzpicture}[font=\footnotesize, edge from parent/.style={draw,thick}]
  \node {$(G,N_{6},1_{N_{6}})$} [grow'=up]
     child {node {$(G_{1},N_{5},\varphi_{1})$}edge from parent[->]
      child {node {$(G_{2},\langle x_{5},x_{6}\rangle,\varphi_{2})$}
       child {node {$(G_{2},G_{2},\varphi_{3})$~~~~~~}}
       child {node {~~~~~~$(G_{2},G_{2},\varphi_{4})$}}
    }};
\end{tikzpicture} \caption{$\mathcal{G}_{N_{6}}$} \end{figure} \noindent The leaves $(G_{2},G_{2},\varphi_{3})$ and $(G_{2},G_{2},\varphi_{4})$ of $\mathcal{G}_{N_{6}}$ correspond to the Shoda pairs $(\langle x_{5},x_{6},x_{3}x_{4}^{2}\rangle, \langle x_{3}x_{4}^{2}x_{6}^{3},x_{5}^{3}\rangle)$ and $(\langle x_{5},x_{6},x_{3}x_{4}^{2}\rangle,\langle x_{5}^3\rangle)$ of $G$ respectively. It turns out that the collection of Shoda pairs corresponding to all the $\mathcal{G}_{N_{i}}$, $1 \leq i \leq 6$, is a complete irredundant set of Shoda pairs of $G$. Futhermore, if $(H,K)$ is any of the Shoda pairs constructed with this process, then, from Theorem \ref{t4}, it follows that $$\footnotesize{\alpha_{(G,H,K)}=\left\{\begin{array}{ll}1/2, & \hbox{if $(H,K)$=$(\langle x_{5},x_{6},x_{3}x_{4}^{2}\rangle, \langle x_{3}x_{4}^{2}x_{6}^{3},x_{5}^{3}\rangle)$ or $(\langle x_{5},x_{6},x_{3}x_{4}^{2}\rangle,\langle x_{5}^3\rangle)$;} \\1, & \hbox{otherwise.}\end{array}\right.}$$

\subsection{Example $2$} Consider the group $G$  generated by $a,~b,~c,~d,~e,~f$ with the following defining relations:
    \begin{center} $a^{2}$=$b^{3}$=$c^{3}$=$d^{3}$=$1$, $a^{-1}ba$=$b^{-1}$, $a^{-1}ca$=$c^{-1}$, $a^{-1}da$=$d$, $b^{-1}cb$=$cd$, $b^{-1}db$=$d$, $c^{-1}dc$=$d$.
                 \end{center}
\noindent There are $8$ normal subgroups of $G$ given as follows: $N_{1}$=$G$, $N_{2}$=$\langle b,c\rangle$, $N_{3}$=$\langle c,d\rangle$, $N_{4}$=$\langle cb,c^{-1}b^{-1}\rangle$, $N_{5}$=$\langle cb^{-1},c^{-1}b\rangle$, $N_{6}$=$\langle b,d\rangle$, $N_{7}$=$\langle d\rangle$ and
 $N_{8}=\langle 1\rangle.$ As before, the directed tree $\mathcal{G}_{N_{1}}$ is just the vertex $(G,G,1_{G})$, which gives Shoda pair $(G, G)$. As $G/N_{2}$ is cyclic, by fact 3, the tree corresponding to $N_{2}$ is as given below:\begin{figure}[H] \centering \begin{tikzpicture}[font=\footnotesize, edge from parent/.style={draw,thick}]
  \node {$(G,N_{2},1_{N_{2}})$} [grow'=up]
     child {node {$(G,G,\varphi)$~~~~}edge from parent[->]}
     ;
\end{tikzpicture} \caption{$\mathcal{G}_{N_{2}}$} \end{figure}
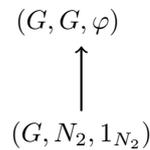 \noindent where $\varphi$ can be taken as any linear character on $G$ with kernel $N_{2}$, which gives Shoda pair  $(G, N_{2})$ of $G$. \par  For $3\leq i\leq 6$, $N_{2}/N_{i}$ is an abelian normal subgroup of maximal order in $G/N_{i}$. Hence we set $\mathcal{A}_{(G, N_{i}, 1_{N_{i}})}= N_{2}$. It turns out that $\mathcal{A}_{(G, N_{i}, 1_{N_{i}})}/N_{i}$ is cyclic and maximal among all the abelian subgroups of $G/N_{i}$. Therefore by fact 2, $\mathcal{G}_{N_{i}}$, $3\leq i\leq 6$, can be described as follows:\\
\begin{figure}[H]
\centering
\begin{minipage}{.25\textwidth}\centering
\begin{tikzpicture}[font=\footnotesize, edge from parent/.style={draw,thick}]
  \node {$(G,N_{3},1_{N_{3}})$} [grow'=up]
     child {node {$(\langle b,c \rangle,\langle b,c \rangle,\varphi_{1})$~~~~}edge from parent[->]}
     ;\end{tikzpicture}
     \captionof{figure}{$\mathcal{G}_{N_{3}}$}
     \end{minipage}%
\begin{minipage}{.25\textwidth}
  \centering
~~\begin{tikzpicture}[font=\footnotesize, edge from parent/.style={draw,thick}]
  \node {$(G,N_{4},1_{N_{4}})$} [grow'=up]
     child {node {$(\langle b,c \rangle,\langle b,c \rangle,\varphi_{2})$~~~~}edge from parent[->]}
     ;
\end{tikzpicture}
\captionof{figure}{$\mathcal{G}_{N_{4}}$}
\end{minipage}%
\begin{minipage}{.25\textwidth}
  \centering
~~~\begin{tikzpicture}[font=\footnotesize, edge from parent/.style={draw,thick}]
  \node {$(G,N_{5},1_{N_{5}})$} [grow'=up]
     child {node {$(\langle b,c \rangle,\langle b,c \rangle,\varphi_{3})$~~~~}edge from parent[->]}
     ;
\end{tikzpicture}
\captionof{figure}{$\mathcal{G}_{N_{5}}$}
\end{minipage}%
\begin{minipage}{.25\textwidth}
  \centering
~~~\begin{tikzpicture}[font=\footnotesize, edge from parent/.style={draw,thick}]
  \node {$(G,N_{6},1_{N_{6}})$} [grow'=up]
     child {node {$(\langle b,c \rangle,\langle b,c \rangle,\varphi_{4})$~~~~}edge from parent[->]}
     ;
\end{tikzpicture}
\captionof{figure}{$\mathcal{G}_{N_{6}}$}
\end{minipage}
\end{figure}
\noindent where $\varphi_{1},~\varphi_{2},~\varphi_{3},~\varphi_{4}$ are any linear characters on $\langle b,c \rangle$ with kernel $N_{3}$, $N_{4}$, $N_{5}$, $N_{6}$ respectively.  The above trees yield Shoda pairs $(\langle b,c \rangle, N_{3})$, $(\langle b,c \rangle, N_{4})$, $(\langle b,c \rangle, N_{5})$, $(\langle b,c \rangle, N_{6})$ of $G$. \par Next, it turns out that $\mathcal{C}l(G,N_{7},1_{N_{7}})$ is empty. Therefore, $\mathcal{G}_{N_{7}}$ is just the vertex $(G,N_{7},1_{N_{7}})$, which does not give any Shoda pair. \par  It now remains to construct the directed tree $\mathcal{G}_{N_{8}}$.  Observe that $N_{3}/N_{8}$ is an abelian normal subgroup of maximal order in $G/N_{8}$. Hence we can set $ \mathcal{A}_{(G,N_{8},1_{N_{8}})}=N_{3}$. In view of fact 5, we have $Cl(G,N_{8},1_{N_{8}})= \{(G_{1}, N_{3}, \varphi)\} $, where $\varphi$ is a linear character on $N_{3}$ with kernel $\langle c\rangle$ and $G_{1}$=$I_{G}(\varphi)$=$\langle a, c,d \rangle.$ As $G_{1}\neq N_{3}$, we further compute $Cl(G_{1}, N_{3}, \varphi)$. Now observe that $G_{1}/\operatorname{ker} \varphi = G_{1}/\langle c\rangle$ is abelian and therefore facts 5 \& 6 yield $$Cl(G_{1}, N_{3}, \psi) = \{(G_{1},G_{1}, \varphi_{1}),~(G_{1},G_{1}, \varphi_{2})\},$$ where $\varphi_{1}$ and $\varphi_{2}$  can be taken to be any linear characters of $G_{1}$ with kernel  $\langle c\rangle$ and  $\langle a,c \rangle$  respectively. The process stops here and the corresponding tree is as follows:
\begin{figure}[H] \centering \begin{tikzpicture}[font=\footnotesize, edge from parent/.style={draw,thick}]
  \node {$(G,N_{8},1_{N_{8}})$} [grow'=up]
     child {node {$(G_{1},\langle c,d \rangle,\varphi)$}edge from parent[->]
       child {node {$(G_{1},G_{1}, \varphi_{1})$~~~~~~}}
       child {node {~~~~~~$(G_{1},G_{1}, \varphi_{2})$}}
    };
\end{tikzpicture}\caption{$\mathcal{G}_{N_{8}}$} \end{figure}
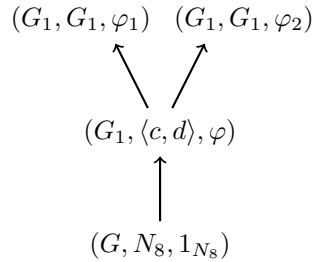 \noindent and it yields Shoda pairs $(\langle a, c,d \rangle, \langle c\rangle)$ and  $(\langle a, c,d \rangle, \langle a,c \rangle)$ of $G$. Again, it turns out that the collection of Shoda pairs constructed with this process is a complete and irredundant set of Shoda pairs of $G$.

\bibliographystyle{amsplain}
\bibliography{Bibliography}

\providecommand{\bysame}{\leavevmode\hbox to3em{\hrulefill}\thinspace}
\providecommand{\MR}{\relax\ifhmode\unskip\space\fi MR }
\providecommand{\MRhref}[2]{%
  \href{http://www.ams.org/mathscinet-getitem?mr=#1}{#2}
}
\providecommand{\href}[2]{#2}
\begin{thebibliography}{10}

\bibitem{BKP}
G.~K. Bakshi, R.~S. Kulkarni, and I.~B.~S. Passi, \emph{The rational group
  algebra of a finite group}, J. Algebra Appl. \textbf{12} (2013), no.~3.

\bibitem{BM}
G.K. Bakshi and S.~Maheshwary, \emph{The rational group algebra of a normally
  monomial group}, J. Pure Appl. Algebra \textbf{218} (2014), no.~9,
  1583--1593.

\bibitem{BM1}
\bysame, \emph{Extremely strong shoda pairs with \texttt{GAP}}, J. Symbolic
  Comput. \textbf{76} (2016), no.~5, 97--106.

\bibitem{wedd}
O.~Broche~Cristo, A.~Herman, A.~Konovalov, A.~Olivieri, G.~Olteanu, {\'A}.~del
  R{\'i}o, and I.~van Geldar, \emph{Wedderga --- wedderburn decomposition of
  group algebras}, Version 4.7.2; (2014), (http://www.cs.st-andrews.ac.uk/
  \~alexk/wedderga).

\bibitem{CR}
C.~W. Curtis and Irving Reiner, \emph{Representation theory of finite groups
  and associative algebras}, Pure and Applied Mathematics, Vol. XI,
  Interscience Publishers, a division of John Wiley \& Sons, New York-London,
  1962.

\bibitem{ND}
Narsingh Deo, \emph{Graph theory with applications to engineering and computer
  science}, Prentice-Hall Series in Automatic Computations, Prentice-Hall,
  Inc., Englewood Cliffs, 1974.

\bibitem{GAP4}
The GAP~Group, \emph{{GAP -- Groups, Algorithms, and Programming, Version
  4.8.6}}, 2016.

\bibitem{HowI}
G.~A. How, \emph{Special classes of monomial groups. {I}}, Chinese J. Math.
  \textbf{12} (1984), no.~2, 115--120.

\bibitem{HowII}
\bysame, \emph{Special classes of monomial groups. {II}}, Chinese J. Math.
  \textbf{12} (1984), no.~2, 121--127.

\bibitem{HowIII}
\bysame, \emph{Special classes of monomial groups. {III}}, Chinese J. Math.
  \textbf{12} (1984), no.~3, 199--211.

\bibitem{BH}
Bertram Huppert, \emph{Character theory of finite groups}, De Gruyter
  Expositions in Mathematics, 25, Walter de Gruyter \& Co., Berlin, 1998.

\bibitem{IM}
I.~Martin Isaacs, \emph{Character theory of finite groups}, Academic Press
  [Harcourt Brace Jovanovich Publishers], New York, 1976, Pure and Applied
  Mathematics, No. 69.

\bibitem{IM1}
I.Martin Isaacs, \emph{Finite group theory}, Graduate Study in Mathematics, 92,
  American Mathematical Society Providence, R I, 2008.

\bibitem{OdRS04}
A.~Olivieri, {\'A}.~del R{\'{\i}}o, and J.~J. Sim{\'o}n, \emph{On monomial
  characters and central idempotents of rational group algebras}, Comm. Algebra
  \textbf{32} (2004), no.~4, 1531--1550.

\bibitem{TY}
T.~Yamada, \emph{The {S}chur subgroup of the {B}rauer group}, Lecture Notes in
  Mathematics, Vol. 397, Springer-Verlag, Berlin, 1974.

\end{thebibliography}

\end{document}